# SYMMETRY VIA SPHERICAL REFLECTION AND SPANNING DROPS IN A WEDGE

JOHN MCCUAN

Abstract. We consider embedded ring-type surfaces (that is, compact, connected, orientable surfaces with two boundary components and Euler-Poincaré characteristic zero) in $\mathbf{R}^3$ of constant mean curvature which meet planes $\Pi_1$ and $\Pi_2$ in constant contact angles $\gamma_1$ and $\gamma_2$ and bound, together with those planes, an open set in $\mathbf{R}^3$. If the planes are parallel, then it is known that any contact angles may be realized by infinitely many such surfaces given explicitly in terms of elliptic integrals. If $\Pi_1$ meets $\Pi_2$ in an angle $\alpha$ and if $\gamma_1 + \gamma_2 > \pi + \alpha$, then portions of spheres provide (explicit) solutions. In the present work it is shown that if $\gamma_1 + \gamma_2 \leq \pi + \alpha$, then the problem admits no solution. The result contrasts with recent work of H.C. Wente who constructed, in the particular case $\gamma_1 = \gamma_2 = \pi/2$, a *self-intersecting* surface spanning a wedge as described above.

Our proof is based on an extension of the Alexandrov planar reflection procedure to a reflection about spheres [7], on the intrinsic geometry of the surface, and on a new maximum principle related to surface geometry. The method should be of interest also in connection with other problems arising in the global differential geometry of surfaces.

## 1. Introduction

We seek to characterize the embedded surfaces of constant mean curvature which span a wedge of opening angle $\alpha$ and which meet the planes $\Pi_1$ and $\Pi_2$ of the wedge in constant contact angles $\gamma_1$ and $\gamma_2$ respectively (see Figure 1). Wente recently constructed an immersed such surface (spanner) for which $\gamma_1 = \gamma_2 = \pi/2$ [17]. Wente's example is a ring-type surface, but self-intersecting.

(A *ring-type surface* is a compact, connected, orientable surface with two boundary components and Euler-Poincaré characteristic zero. Ring-type surfaces are also known as topologically annular.)

In the present work we establish the following non-existence result.

**Theorem 1.** *If $\gamma_1 + \gamma_2 \leq \pi + \alpha$, then there are no embedded ring-type surfaces of constant mean curvature spanning a wedge of angle $\alpha$ and maintaining constant contact angles, $\gamma_1$ and $\gamma_2$, with the wedge on each component of their boundary.*

This work was supported in part by NASA grant NAS8-39225 (Gravity Probe B Relativity Mission). Research at MSRI supported in part by NSF grant DMS 9022140.

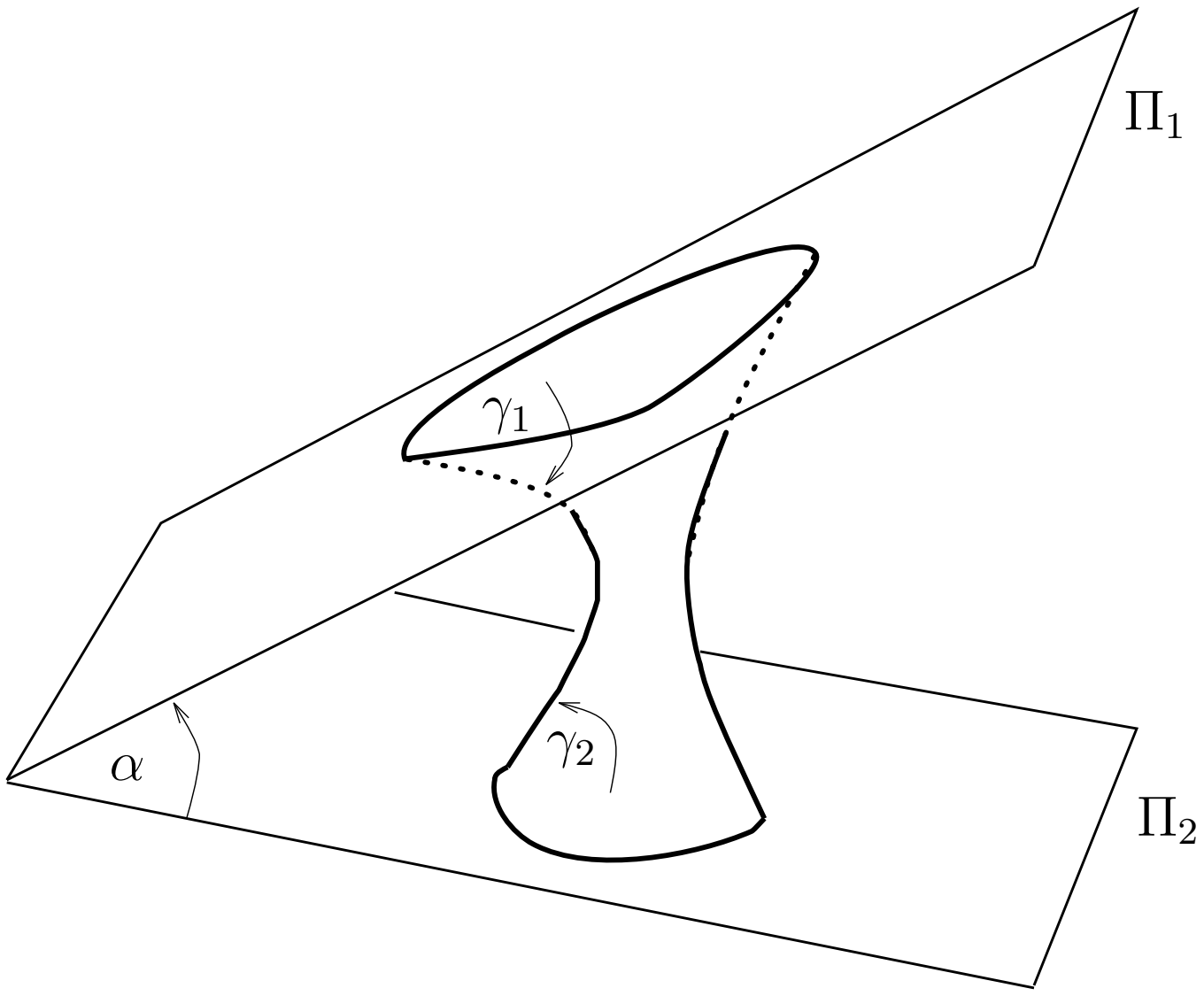

FIGURE 1. Spanning Drop

The contact angle condition of Theorem 1 cannot be improved, as spanners can be found explicitly using spheres in the complementary case:

**Theorem 2.** *A spherical spanner exists if and only if* $\gamma_1 + \gamma_2 > \pi + \alpha$. *Letting* $A > 0$ *denote either enclosed volume or mean curvature, the family of all spanning spheres may be indexed uniquely by the set of 4-tuples* $(\gamma_1, \gamma_2, \alpha, A)$ *for which this existence criterion is satisfied.* (See Figure 2.)

Whether or not non-spherical embedded spanners exist satisfying $\gamma_1 + \gamma_2 > \pi + \alpha$ remains an open question.

Any embedded spanning surface (spherical or otherwise) corresponds to an equilibrium for the free surface of a volume of liquid in the absence of gravity. A considerable literature has developed treating the stability of equilibrium spanning surfaces between parallel planes; these are the so-called liquid bridges (see for example [2, 14, 4, 12, 18]). It is a remarkable corollary of Theorem 1 that, although there are many equilibrium spanners (even stable ones) between parallel planes satisfying $\gamma_1 + \gamma_2 \leq \pi$, if the planes make *any* angle $\alpha \neq 0$ with each other, then not only stability but the existence of equilibria ceases.

Experimentally, one may observe that fluid drops with $\gamma < \pi/2$ or soap films ($\gamma = \pi/2$) that span a wedge are unstable and tend to the vertex. In fact, these investigations were prompted by the common practice (in the engineering of liquid propellant tanks) of placing a wedge shaped partition with its vertex on the fluid outlet to ensure a supply of fluid there. The presence of fluid at the vertex is said to arise from the wedge's ability to "force bubbles out." Our result is a first step in mathematical verification of this phenomenon (see Figure 3).

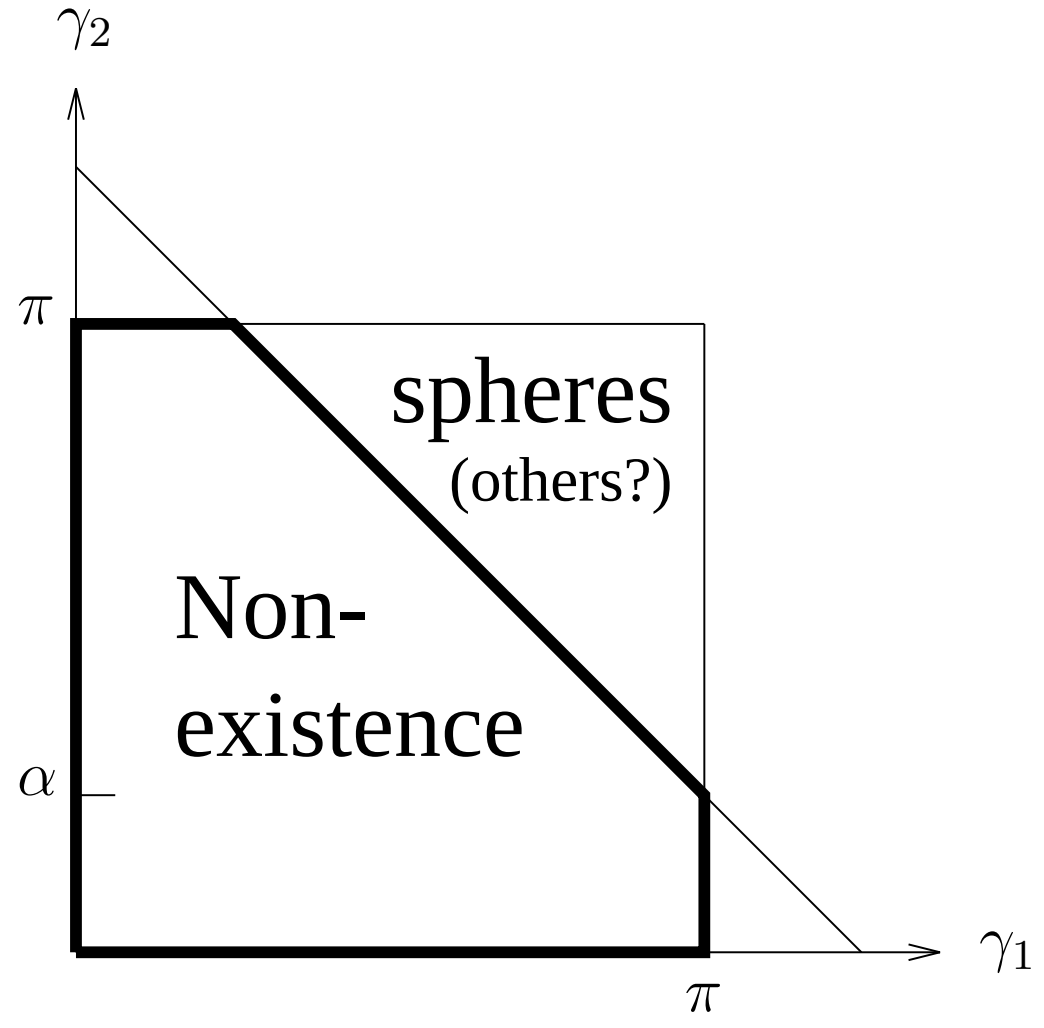

FIGURE 2. Existence and Non-existence

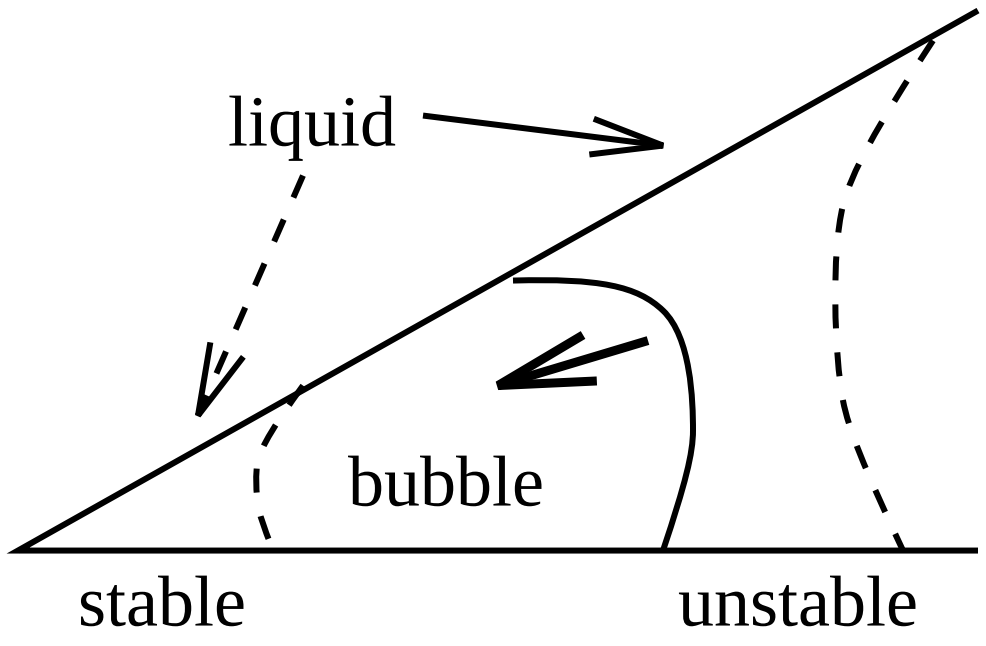

FIGURE 3. "forcing bubbles out"

The proof of the main result depends on an extension of the reflection method of A.D. Alexandrov to non-planar reflecting surfaces, namely spheres [7]. The main results in this connection are stated in §2. Furthermore, it is observed that the *maximum principle of spherical reflection* [7] essentially reduces the problem to consideration of the boundary behavior of the surface. This observation is applied in detail in §3, and methods of H. Hopf and, hence, the topological assumptions are employed to complete the proof of Theorem 1.

It will be observed that previous applications of reflection, cf. [1, 11, 15, 5, 10], do not rely on topological assumptions as does Theorem 1. In principle, this is true for spherical reflection as well, and therefore it is natural to seek an extension of Theorem 1 to surfaces of higher topology. In a separate paper [8], using different methods, the following partial extension in this direction is given (along with some other extensions of less immediate physical interest).

**Theorem 3.** *If $\gamma_1, \gamma_2 \leq \pi/2$, then there are no embedded constant mean curvature spanners (in any wedge).*

Notice, however, that the sharp contact angle condition has not been obtained.

I would like to acknowledge the guidance and encouragement of my advisor Robert Finn during the course of this work. I am also thankful for helpful conversations with Rick Schoen, Leon Simon and Brian White.

## 2. Spherical Reflection of Spanning Drops

Let $\mathcal{S}$ be the free surface of a volume of fluid which adheres to the planes of a wedge. According to our assumptions on such a surface, which are discussed in detail in Appendix A, there is an open set $\mathcal{D}$ in $\mathbf{R}^3$ whose boundary is composed of $\mathcal{S}$ and the regions of adherence $B$ which are topological discs in the planes $\Pi$ of the wedge; see also Figure 1. We begin the reflection procedure by fixing an origin $0 \in \mathbf{R^3}$ at a point on the vertex of the wedge. Next, consider a large sphere $S_\rho = \partial B_\rho(0)$ such that $\bar{\mathcal{D}} \subset B_\rho(0) \equiv \{\xi \in \mathbf{R}^3 : |\xi| < \rho\}$. As we decrease the radius $\rho$ we come to the first sphere $S_{\rho_0}$ which intersects $\mathcal{S}$. As we decrease $\rho$ below $\rho_0$ we apply the (reflection) map

$$\psi : X \longmapsto \frac{\rho^2}{|X|^2} X \tag{1}$$

to the portion of $\mathcal{S}$ through which $S_\rho$ has passed. To be precise, let $\mathcal{S}^- = \mathcal{S}^-(\rho) \equiv \{X \in \mathcal{S} : |X| \geq \rho\}$; the reflection of $\mathcal{S}^-$ is $\hat{\mathcal{S}} = \hat{\mathcal{S}}(\rho) \equiv \psi \mathcal{S}^-$. A description of the of the spatial disposition of $\mathcal{S}^+(\rho) \equiv \{X \in \mathcal{S} : |X| \leq \rho\}$ and $\hat{\mathcal{S}}$ during the reflection procedure is as follows.

**Theorem 4.** *There is some radius $\rho_1$ with $0 < \rho_1 < \rho_0$ such that for each $\rho \in (\rho_1, \rho_0)$ every point $X$ in $\mathcal{S}^-(\rho)$ satisfies exactly one of the following.*

**NT1.:** $|X| > \rho$, $X \in \text{int } \mathcal{S}$, *and* $(1-\delta)X \in \mathcal{D}$ *for each* $\delta \in (0, 1 - \rho^2/|X|^2]$.[1]

**NT2.:** $|X| > \rho$, $X \in \partial\mathcal{S}$, *and* $(1-\delta)X \in B$ *for each* $\delta \in (0, 1 - \rho^2/|X|^2]$.

**NT3.:** $|X| = \rho$, $X \in \text{int } \mathcal{S}$ *and* $X \cdot N < 0$ *where* $N$ *is the normal pointing into* $\mathcal{D}$.

**NT4.:** $|X| = \rho$, $X \in \partial\mathcal{S}$, *and* $X \cdot n < 0$, *where* $n$ *is the inward normal to* $\partial\mathcal{S}$ *in* $\Pi$.

*For $\rho = \rho_1$, however, at least one of the following holds for some $X \in \mathcal{S}^-$.*

[1] Notice that the right endpoint in this interval corresponds to $\hat{X}$, the reflection of $X$ about $S_\rho$.

**T1.:** *$|X| > \rho$, $X \in \text{int}\,\mathcal{S}$, and $\hat{\mathcal{S}}$ is tangent to $\mathcal{S}^+$ at $\hat{X}$.*
**T2.:** *$|X| > \rho$, $X \in \partial\mathcal{S}$, and $\hat{\mathcal{S}}$ is tangent to $\mathcal{S}^+$ at $\hat{X}$.*
**T3.:** *$|X| = \rho$, $X \in \text{int}\,\mathcal{S}$ and $N$ is tangent to $S_{\rho_1}$ at $\hat{X} = X$.*
**T4.:** *$|X| = \rho$, $X \in \partial\mathcal{S}$, and $n$ is tangent to $S_{\rho_1}$ at $\hat{X} = X$.*

PROOF: See Appendix A.

As discussed in Appendix A, if the mean curvature of $\hat{\mathcal{S}}$ with respect to the reflection $\hat{N}$ of the normal $N$ is bounded above by the mean curvature $H$ of $\mathcal{S}$, then a touching principle leads to the conclusion $\hat{\mathcal{S}}(\rho_1) \equiv \mathcal{S}^+(\rho_1)$. Our discussion, therefore, centers on the following global property of spanning surfaces.

**Definition 5** (Boundedness Property)**.** *The CMC*[2] *surface $\mathcal{S}$ is said to have the* Boundedness Property *if for each $\rho \geq \rho_1$ and each $X \in \mathcal{S}$ with $|X| \geq \rho$ there holds the inequality*

$$\hat{H}(X, \rho) \leq H$$

*where $\hat{H}(X, \rho)$ is the mean curvature of $\hat{\mathcal{S}}$ at the image of $X$ under the map (1).*

In the proof of Theorem 4 we obtain the following additional information.

**Lemma 6.** *Consider $X \in \mathcal{S}^-(\rho_1)$, i.e. $X \in \mathcal{S}$ and $|X| \geq \rho_1$.*
**(i):** *$X \cdot N \leq 0$.*
**(ii):** *Among $X \in \text{int}\,\mathcal{S}$, equality can hold only for $|X| = \rho_1$.*
**(iii):** *Among $X \in \partial\mathcal{S}$ with $0 < \gamma_1, \gamma_2 < \pi$, equality can hold only for $|X| = \rho_1$.*

By combining part (i) of this result with the formula

$$\hat{H}(X, \rho) = \frac{1}{\rho^2}(|X|^2 H + 2X \cdot N) \tag{2}$$

(demonstrated in [7]) we have

**Theorem 7.** *Any spanning surface of non-positive constant mean curvature satisfies the boundedness property.*

Thus, we turn our attention to surfaces with positive mean curvature. A fundamental observation in [7] is the following.

**Theorem 8.** *If $H > 0$ and $\rho \in [\rho_1, \rho_0)$ is fixed, then $\Delta \hat{H}(X, \rho) > 0$ where $\Delta$ is the intrinsic Laplacian on $X \in \mathcal{S}^-(\rho)$.*

Thus, $\hat{H}$ is subharmonic and achieves its maximum at a point $\bar{X} = \bar{X}(\rho) \in \partial\mathcal{S}^-$. If $|\bar{X}| = \rho$, then according to (2) and Lemma 6 for each $X \in \mathcal{S}^-(\rho)$,

[2] Constant Mean Curvature

$$\begin{aligned}\hat{H}(X,\rho) &\le \hat{H}(\bar{X},\rho)\\ &= H + \frac{2}{\rho^2}\bar{X}\cdot\bar{N}\\ &\le H.\end{aligned}$$

This is the condition required by the boundedness property. In this way, it is only necessary to consider the case when $|\bar{X}| > \rho$ for which $\bar{X} \in \partial\mathcal{S}$. It is this observation that makes the topological assumption useful.

## 3. Ring Type Spanners

A ring type surface of constant mean curvature which meets an umbilic surface (plane or sphere) at a constant contact angle can be parameterized by a single conformal curvature coordinate $X$. This is shown in certain cases in [16] and [9]. We give a derivation of this fact in general which is due to R. Finn.

Let the surface in question be $\mathcal{S}$ and its orientation $N$. Let $Y$ be a conformal representation of $\mathcal{S}$ on a unit annulus $A$ in the $z = x + iy-$ plane (see Figure 4).[3] Let $w = \log z$ map $A$ to a rectangle $\Omega$ in the $w = u + iv-$ plane and notice that the

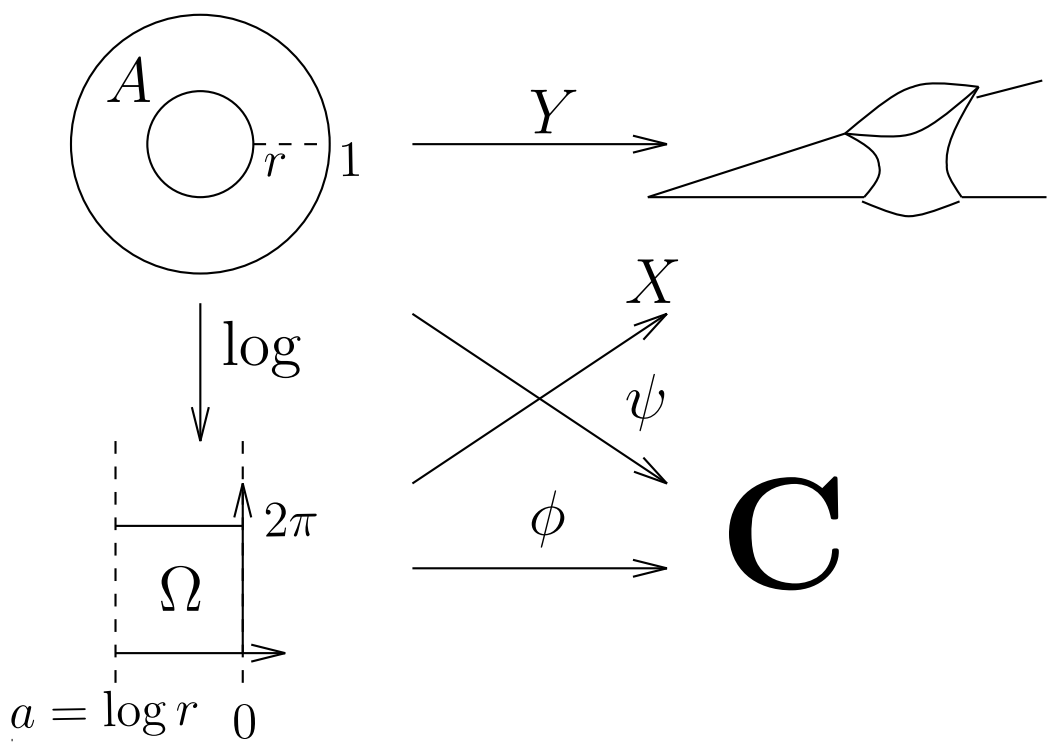

Figure 4. Conformal Representation of Ring Type Surfaces

inverse (exponential) is periodic in the strip $\log r < u < 0$ which we will also refer to as $\Omega$. The map $X = Y \circ \exp : \Omega \to \mathcal{S}$ is again a conformal representation of $\mathcal{S}$. By the *Theorem of Joachimstahl* [3], however, the boundary curves on $\mathcal{S}$ are curvature lines for $\mathcal{S}$. Thus, the coefficient $f = X_{uv}\cdot N$ of the second fundamental form vanishes on the boundary of the strip. On the other hand, $\phi \circ \log$ defines an analytic function on $A$ where $\phi \equiv (e-g) - 2if,\ e = X_{uu}\cdot N,$ and $g = X_{vv}\cdot N$ [6, pg. 139 Lemma 2.2]. $\mathrm{Im}(\phi\circ\log)$ is, therefore, harmonic and vanishes on $\partial A$. Thus, $\mathrm{Im}(\phi\circ\log) \equiv 0$ and

[3] We assume that $Y$ is smooth on $\bar{A}$, whence $X$ is smooth on $\bar{\Omega}$.

$\mathrm{Re}(\phi \circ \log) \equiv c$ (constant). Thus, $f \equiv 0$ so that we have curvature coordinates on $\Omega$, and moreover $e - g \equiv c$ (constant).

In summary, the coefficients of the first and second fundamental forms satisfy

$$E = X_u \cdot X_u = X_v \cdot X_v = G; \;\; F = X_u \cdot X_v = 0,$$

and

$$e - g = X_{uu} \cdot N - X_{vv} \cdot N = c \quad \text{(constant)}; \;\; f = X_{uv} \cdot N = 0.$$

Therefore, $2H = (e+g)/E$, and we see that

$$e = EH + \frac{c}{2}; \;\; g = EH - \frac{c}{2}.$$

Moreover, $k_1 + k_2 = 2H$ and $k_1 k_2 = K = eg/E^2$, so that the principal curvatures are given by

$$k_1 = \frac{e}{E} = H + \frac{c}{2E}; \;\; k_2 = \frac{g}{E} = H - \frac{c}{2E}.$$

Finally, if we take the inward normal $N$ along the coordinate lines (which are curvature lines), then

$$N_u = -k_1 X_u; \;\; N_v = -k_2 X_v.$$

Having noted these auxiliary facts, we can state the two main results of this section.

**Lemma 9.** *The boundedness property holds for a spanning surface of ring type unless $H > 0$ and $c > 0$.*

**Lemma 10.** *If $H, c > 0$, then $\gamma_1 + \gamma_2 > \pi + \alpha$.*

We see from Lemma 9 that, in particular, the boundedness property is intrinsic for ring type spanners in the sense that it is independent of the particular origin of reflection along the vertex of the wedge. Thus, if $\mathcal{S}$ satisfies the boundedness property, then centered at each point $x$ along the vertex of the wedge there is a sphere $S_\rho(x)$ about which $\mathcal{S}$ is symmetric, i.e., invariant under reflection. It is shown in [7] that such a surface is a portion of a sphere as described in Theorem 2. Combining this fact with Lemma 10 evidently yields Theorem 1.

PROOF of Theorem 9: According to Theorem 7 the content of this lemma is that if $\mathcal{S}$ fails to satisfy the boundedness property, then $c > 0$.[4] This is easily verified if one of the contact angles is 0 or $\pi$.

**Lemma 11.** *If the contact angle $\gamma$ is 0 or $\pi$ on the boundary curve $\beta$, then*

$$H = \frac{c}{2E}$$

*on $\beta$. In particular, if $H > 0$ (which it must be for the boundedness property to fail), then $c = 2EH > 0$.*

[4] Henceforth, we assume $H > 0$.

PROOF: Recall that by Joachimstahl's Theorem, the boundary curve $\beta$ is a curvature line with *normal* curvature

$$k_2 = H - \frac{c}{2E}.$$

For these extremal contact angles, however, the normal $N$ is perpendicular to $\Pi$, and since $\beta \subset \Pi$, its normal curvature is 0.

For the remaining contact angles ($0 < \gamma_1, \gamma_2 < \pi$), some care is required in applying the results of §2. Let

$$\Omega_1 \equiv \{(u, v) \in \Omega : |X(u, v)| > \rho_1\}.$$

We have observed that $\hat{H}$ takes its supremum on $\Omega_1$ at a point $(\bar{u}, \bar{v}) \in \partial\Omega_1$.

**Corollary 12.** *Let $\bar{X} \equiv X(\bar{u}, \bar{v})$, then for any $\rho > 0$ and $|X| > \rho_1$*

$$\hat{H}(X, \rho) \leq \hat{H}(\bar{X}, \rho).$$

PROOF: According to (2)

$$\frac{1}{\rho_1^2}(|X|^2 H + 2X \cdot N) \leq \frac{1}{\rho_1^2}(|\bar{X}|^2 H + 2\bar{X} \cdot \bar{N}).$$

Multiply both sides by $\rho_1^2/\rho^2$.

At this point we postulate the failure of the boundedness property and derive a contradiction in the case $c < 0$. The next two results give some important information about how this failure must come about.

**Lemma 13.** *There is some $\bar{\rho}$ and some $\bar{X} \in \mathcal{S}^-(\bar{\rho})$ with $\rho_1 < \bar{\rho} < |\bar{X}|$ such that*

(i): *$\hat{H}(\bar{X}, \bar{\rho}) = H$.*
(ii): *$\hat{H}(\bar{X}, \rho) > H$ if $\rho_1 \leq \rho < \bar{\rho}$.*
(iii): *$\hat{H}(X, \rho) \leq H$ if $|X| \geq \rho \geq \bar{\rho}$.*

NOTE: $(\bar{X}, \bar{\rho})$ gives, in some sense, the "first time" that the boundedness property fails.

PROOF of Lemma 13: $\bar{X}$ is defined in Corollary 12. For each $\rho \geq \rho_1$, $\{X \in \mathcal{S} : |X| \geq \rho\}$ may have points with $|X| = \rho$ and others with $|X| > \rho \geq \rho_1$. As for the first points,

$$\hat{H}(X, \rho) = H + \frac{2}{\rho^2} X \cdot N \leq H$$

by Lemma 6. The second points satisfy $\hat{H}(X, \rho) \leq \hat{H}(\bar{X}, \rho)$ by Corollary 12. We conclude that for some $\rho \geq \rho_1$, $\hat{H}(\bar{X}, \rho) > H > 0$. In particular,

$$|\bar{X}|^2 H + 2\bar{X} \cdot \bar{N} > 0.$$

Therefore, the equation $\hat{H}(\bar{X},\rho) = H$ has some (unique) solution $\rho = \bar{\rho}$. In fact,

$$\bar{\rho} = \sqrt{\frac{1}{H}(|\bar{X}|^2 H + 2\bar{X}\cdot\bar{N})}.$$

The equality (i) is clearly satisfied. Moreover, by the evident monotonicity of $\hat{H}(\bar{X},\rho)$ in $\rho$, we conclude that $\bar{\rho} > \rho_1$ and that (ii) holds.

Since $\bar{\rho} > \rho_1$, Corollary 12 applies to all $(X,\rho)$ with $|X| \geq \rho \geq \bar{\rho}$. Hence,

$$\begin{aligned}
\hat{H}(X,\rho) &\leq \hat{H}(\bar{X},\rho) \\
&= \frac{1}{\rho^2}(|\bar{X}|^2 H + 2\bar{X}\cdot\bar{N}) \\
&\leq \frac{1}{\bar{\rho}^2}(|\bar{X}|^2 H + 2\bar{X}\cdot\bar{N}) \\
&= \hat{H}(\bar{X},\bar{\rho}) \\
&= H.
\end{aligned}$$

This is condition (iii).

Finally, since $|\bar{X}| \geq \rho_1$, $\bar{X}\cdot\bar{N} \leq 0$, and it cannot be the case that $|\bar{X}| < \bar{\rho}$ (just write out condition (i)). Therefore, $|\bar{X}| > \rho_1$, and the equality condition of Lemma 6 applies. This means that the condition $|\bar{X}| = \bar{\rho}$ implies

$$\begin{aligned}
\hat{H}(\bar{X},\bar{\rho}) &= H + \frac{2}{\bar{\rho}^2}\bar{X}\cdot\bar{N} \\
&< H
\end{aligned}$$

again contradicting (i). This establishes that $|\bar{X}| > \bar{\rho}$.

The following observation is crucial for us.

**Lemma 14.** $\bar{X} \in \partial\mathcal{S}$.

PROOF: We know that $(\bar{u},\bar{v}) \in \partial\Omega_1$. The image if $\partial\Omega_1$ under $X$ consists of points in $\partial\mathcal{S}$ and points on $S_{\rho_1}$. $\bar{X}$ is not in the latter set by Lemma 13.

Let us summarize our position. We are proving Theorem 9 which means verifying that the boundedness property holds under "most" circumstances. We have reduced our consideration to the contact angles $0 < \gamma_1, \gamma_2 < \pi$ and shown that the boundedness property holds unless $H > 0$ and $c \neq 0$. We now rule out the possibility that $c < 0$. Recall that we have postulated a violation of the boundedness property which occurs for the "first time" at $(\bar{X},\bar{\rho})$ as described by Lemma 13. Thus, we obtain a contradiction via

**Theorem 15.** *If $c < 0$, then there is some point $\bar{\bar{X}} \in \partial\mathcal{S}$ with $|\bar{\bar{X}}| \geq \bar{\rho}$ such that $\bar{\bar{X}}$ is a* (T4) *touching point.*

PROOF: The idea of the proof is simple.

**(i):** *According to the reflected mean curvature formula (2), the osculating sphere (the sphere of radius* $1/H$ *with center at* $\mathcal{C}(X)$*) reflects to a sphere of radius* $|1/\hat{H}|$*. The osculating sphere at* $\bar{X}$*, therefore, is invariant under reflection about* $S_{\bar{\rho}}$*.*

**(ii):** *The osculating sphere intersects* $\Pi$ *(the plane containing* $\bar{X}$*) in a circle,* $C$*, which is invariant to the reflection.*

**(iii):** *The component of* $\partial\mathcal{S}$ *containing* $\bar{X}$ *is a curve* $\beta$ *that curves "more tightly" than* $C$*. By comparing* $\beta$ *and* $C$ *we can obtain a* (T4) *touching point.*

Some of the details:

Let us first note two consequences of the contact angle condition. As noted above $\beta$ has *normal* curvature $k_2$. Therefore, $\beta$ has curvature $k$ given by $k_2 = k \sin\gamma$ as a plane curve in $\Pi$ with respect to its inward normal $n$. Also, the circle

$$C \equiv S_{\frac{1}{H}}(\mathcal{C}(\bar{X})) \cap \Pi$$

has radius $(1/H)\sin\gamma$ and curvature

$$\tilde{k} \equiv \frac{H}{\sin\gamma}.$$

Now, since $k_2 = H - c/2E > H$, we have the relation

$$k > \tilde{k}. \tag{3}$$

As noted, since $\hat{H}(\bar{X}, \bar{\rho}) = H$, reflection about the circle $\partial B_{\bar{\rho}}(0)$ in the plane $\Pi$ leaves $C$ invariant. For this to be true, it is necessary that 0 be outside of $C$ and that the two unique tangent segments from 0 to $C$ be of length $\bar{\rho}$ (see Figure 5). The points of tangency of these two tangent segments divide $C$ into two arcs, and since $|\bar{X}| > \bar{\rho}$, $\bar{X}$ is in the outer one. Without loss of generality, we may take the $y$-axis through the center of $C$ and $\bar{X}$ in the first quadrant. Moreover, we assume $\beta$ is a parameterization by arclength and that $\tilde{\beta}$ is a parameterization of $C$ by arclength so that

$$\beta(0) = \tilde{\beta}(0) = \bar{X}$$

and

$$\bar{X} \cdot \dot{\beta}(0) = \bar{X} \cdot \dot{\tilde{\beta}}(0) \leq 0$$

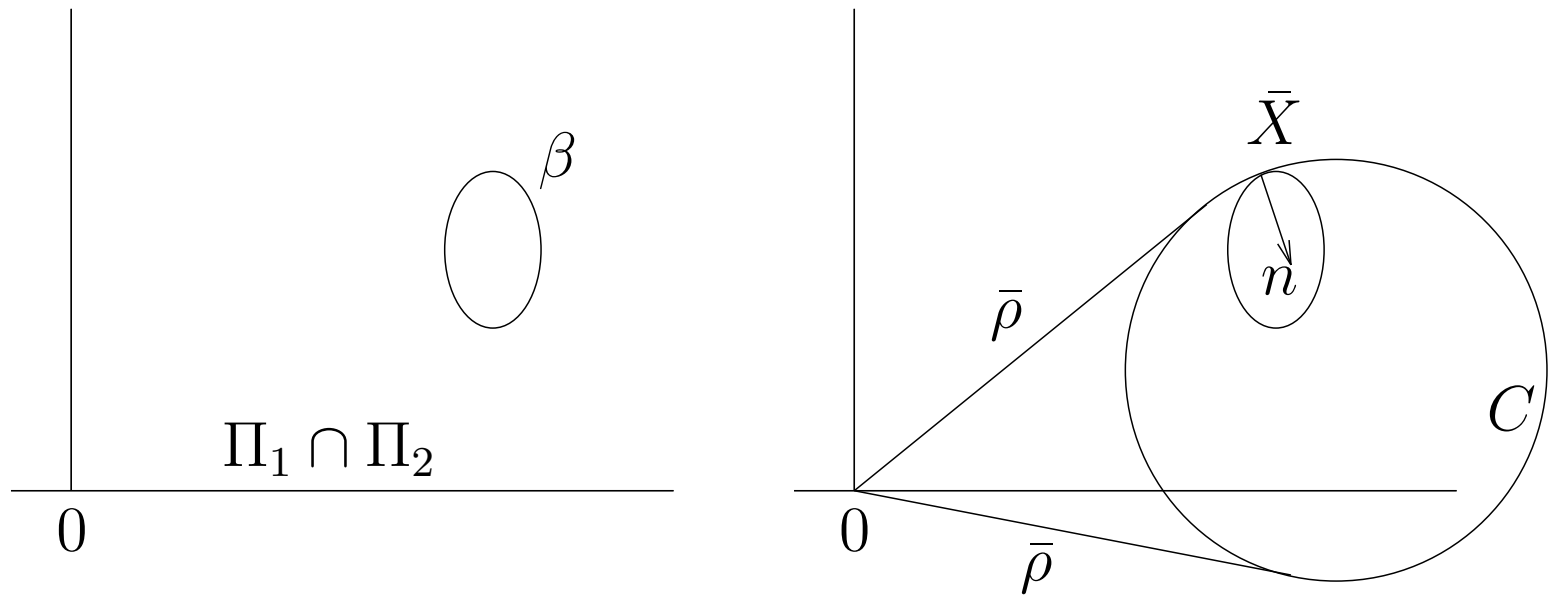

Figure 5. The Plane $\Pi$ (first view)

(see Figure 6). We are interested in the quantities $\phi = |\beta|^2$ and $\psi = (1/k)\beta \cdot \ddot{\beta}$. The following corollary of Lemma 6 explains why.

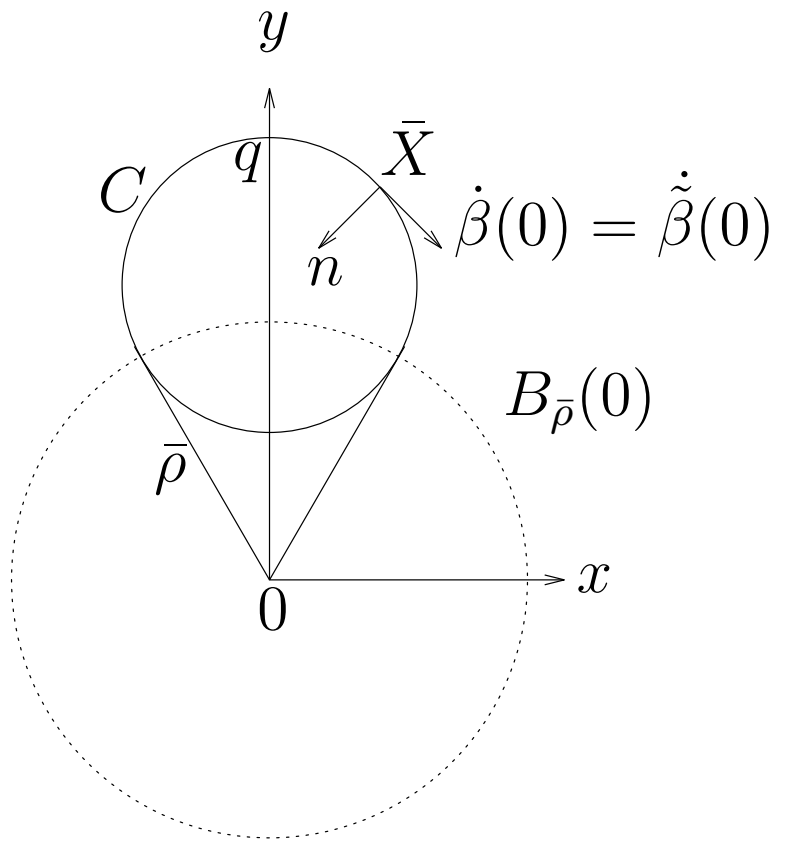

FIGURE 6. The Plane Π (second view)

**Corollary 16.** *If* $0 < \gamma < \pi$, $X \in \partial\mathcal{S}$, *and* $X \cdot N = 0$, *then* $X \cdot n = 0$.

In particular, vanishing of $\psi$ indicates a (T4) touching point. Our contradiction arises from the fact that $\psi$ vanishes at a point where $\phi > \rho_1^2$. In fact, $\psi(0) < 0$, and it can be shown that, at least locally, $\psi$ increases and $\phi$ decreases. On the other hand, by comparison to $C$ using the inequality (3) and the explicit expression for $\beta$ in terms of an integral involving its curvature, one sees that these monotonicity properties persist until $\psi$ vanishes. Furthermore, at this point $\phi \geq \bar{\rho}^2 > \rho_1^2$.

Thus, we have shown that $H, c > 0$. It is a general consequence of these inequalities that $\gamma_1 + \gamma_2 < \pi + \alpha$. In order to show this we turn to planar reflection.

PROOF of Lemma 10: A CMC spanning surface possesses a plane of symmetry $\Sigma$ with normal parallel to the vertex, $l$, of the wedge. This plane is obtained via the conventional reflection procedure and, as a consequence, divides $\mathcal{S}$ into two halves each of which is a graph over $\Sigma$. Notice that the wedge intersects $\Sigma$ in two half lines $\Pi_1$ and $\Pi_2$ which emanate from the point of intersection $l_0$ of $l$ with $\Sigma$. In each of these half lines $\Pi$, there is an open segment $B$ which is the intersection of the wetted region $B$ with $\Sigma$. The endpoints of these intervals $a$ and $b$ (which may be a single point in the case $\gamma = \pi$) are the intersection of $\partial\mathcal{S}$ with $\Sigma$.

Let us fix a frame of reference $x_1, x_2, x_3$ with the $x_1$-axis along $l$, $B$ lying along the positive $x_2$-axis, and the inward normal, $N'$, to $B$ along the positive $x_3$-axis (see Figure 7). Since the normal $N$ to $\mathcal{S}$ along $\mathcal{S} \cap \Sigma$ is invariant under reflection in

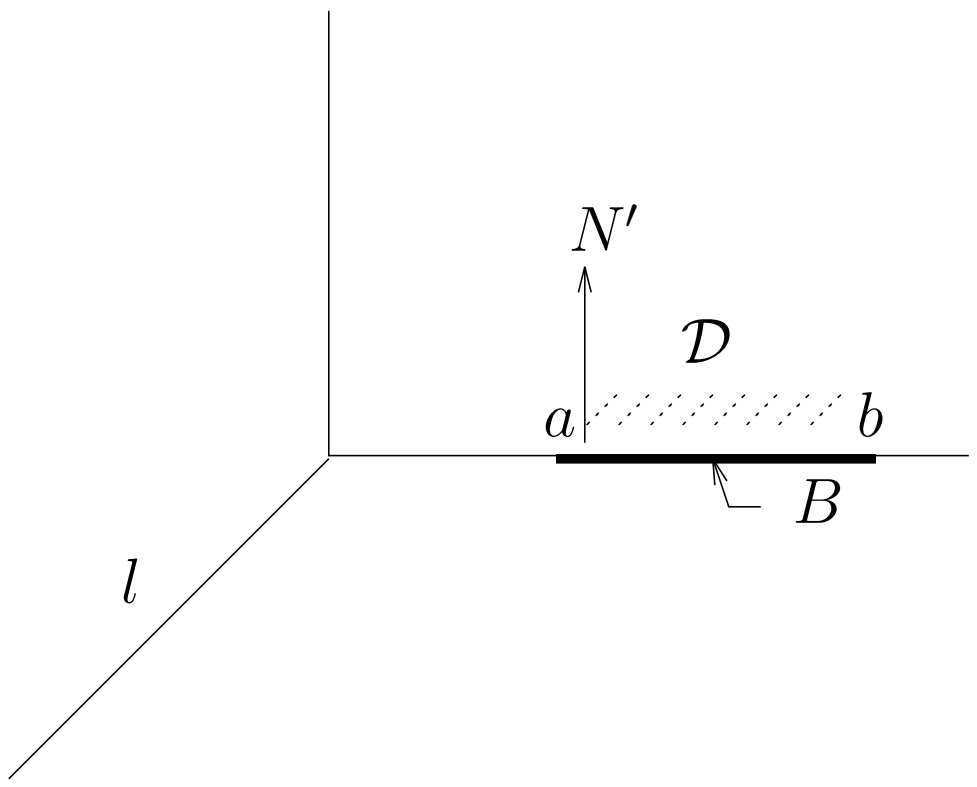

FIGURE 7. The Plane of Symmetry $\Sigma$

$\Sigma$, it must be the case that $N$ lies in $\Sigma$. $\mathcal{S}$ is therefore transverse to $\Sigma$ (at $a$ for example), and intersects $\Sigma$ in a unique smooth curve $c_a$ with endpoint $a$. For points $X \in$ int $\mathcal{S} \cap \Sigma$ the corresponding curve through $X$ contains $X$ in its interior. Thus, the curve, $c_a$, starting at $a$ continues until it reaches some other endpoint $a'$ or $b'$. In this way, the endpoints $\partial\mathcal{S} \cap \Sigma$ are connected to each other in pairs by two smooth curves (see Figure 8). Since $\mathcal{S}$ meets $\Sigma$ orthogonally along $c_a$ and $c_b$ these curves are

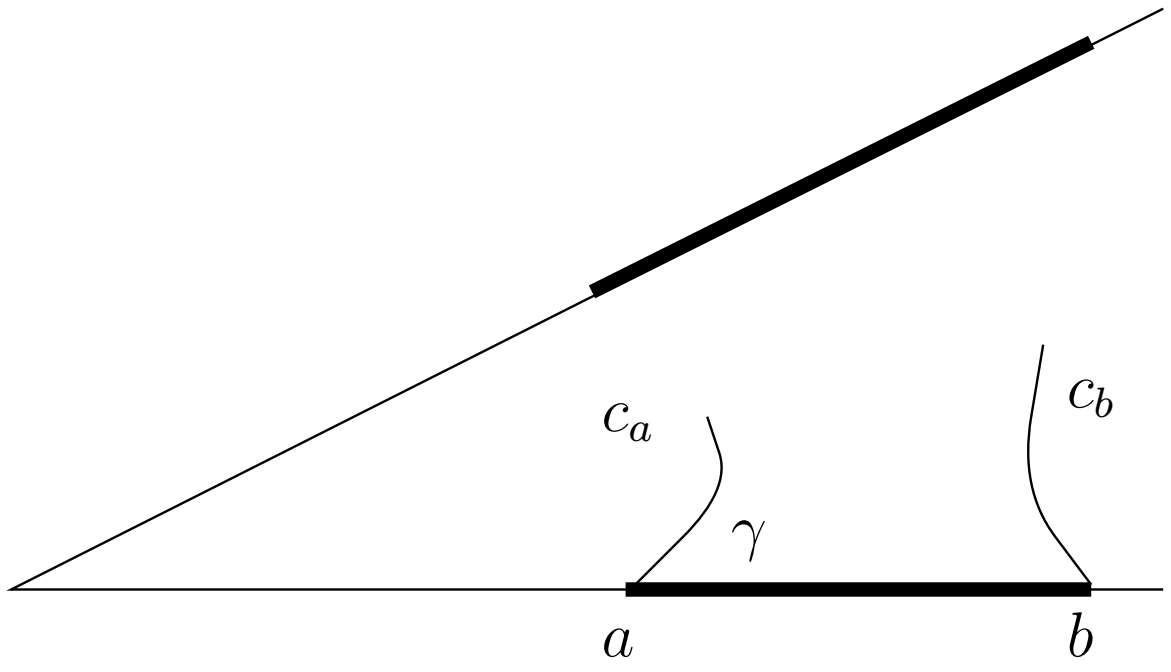

FIGURE 8. The Plane $\Sigma$

curvature lines which have normal curvatures (and hence curvatures) given by

$$k_1 = H + \frac{c}{2E} > H.$$

It follows from this that $\gamma \neq 0$. Moreover, if we let

$$p^\theta \equiv \{(0, \tau \cos\theta, \tau \sin\theta) : \tau > 0\},$$

then for small, positive $\theta$, $p^\theta \cap \{c_a\}$ consists of a unique point. In this way, $c_a$ and $c_b$ can be locally parameterized by $\theta$. Also, for $\theta$ small and positive the straight line segment between $c_a(\theta)$ and $c_b(\theta)$ will lie entirely in $\mathcal{D}$. We define $\theta_0$ to be the largest angle for which this condition holds (and the parameterizations remain valid).

An analysis of the (plane) geometry of these curves as described above results in the condition $\gamma_1 + \gamma_2 < \pi + \alpha$. For example, if $c_a(\theta_0)$ is in $\partial\mathcal{S}$, then draw the chord from $a$ to $c_a(\theta_0)$, and consider the triangle $\triangle\, 0\, c_a(\theta_0)\, a$. Due to the curvature of $c_a$, the contact angle at $a$, $\gamma = \gamma_1$, is greater than the exterior angle at $a$. Thus, if $\phi$ is the angle of the triangle at $c_a(\theta_0)$, then

$$\begin{aligned} \gamma_1 &> \theta_0 + \phi \\ &= \alpha + \phi. \end{aligned}$$

On the other hand, $\gamma_2 > \pi - \phi$. Therefore,

$$\gamma_1 > \alpha + \pi - \gamma_2.$$

This concludes the proof of Lemma 10 and, hence, that of Theorem 1.

## Appendix A. Some Touchy Points

In this section we present the modifications of standard reflection techniques necessary for spherical reflection of spanning drops. The discussion of §2 applies to compact spanning surfaces with arbitrary genus and any number of boundary components. Our discussion below treats some of these surfaces.

A spanning drop, for our purposes, spans a *wedge*. For convenience we fix the vertex of our wedge to be the $X_3$-axis, and *define* a wedge to be a family of open half planes $\Pi_1, \Pi_2, \ldots, \Pi_k$ whose boundaries are all the $X_3$-axis.

We require the drop surface $\mathcal{S}$ to be a compact, connected, orientable, CMC surface with $k$ boundary components $\beta_1, \beta_2, \ldots, \beta_k$ each lying entirely in the open half plane of the same index. Notice that we did not require the half planes to be distinct. We also explicitly assume that $\partial\mathcal{S}$ is smooth although this can be proved, at least in some cases, under a reasonable stability assumption; see [13] and references therein. Each $\beta_j$, furthermore, should be a simple closed curve along which the drop *adheres* to the plane $\Pi_j$. In order to understand what we mean by adherence, let $B_j$ be the topological disk in $\Pi_j$ bounded by $\beta_j$. We require that three conditions hold. First,

A1. $\mathcal{S} \cup B_1 \cup B_2 \cup \cdots \cup B_k$ is a piecewise smooth simple closed surface all of whose singularity lies in $\cup\{\beta_j\}$ (see Figure 9).

Note that by the theorem of Jordan and Brouwer, [6, pg. 100], $\cup B_j \cup \mathcal{S} = \partial\mathcal{D}$ for some open set $\mathcal{D} \subset \mathbf{R}^3$. We will call $\mathcal{D}$ the drop region. Secondly, we want

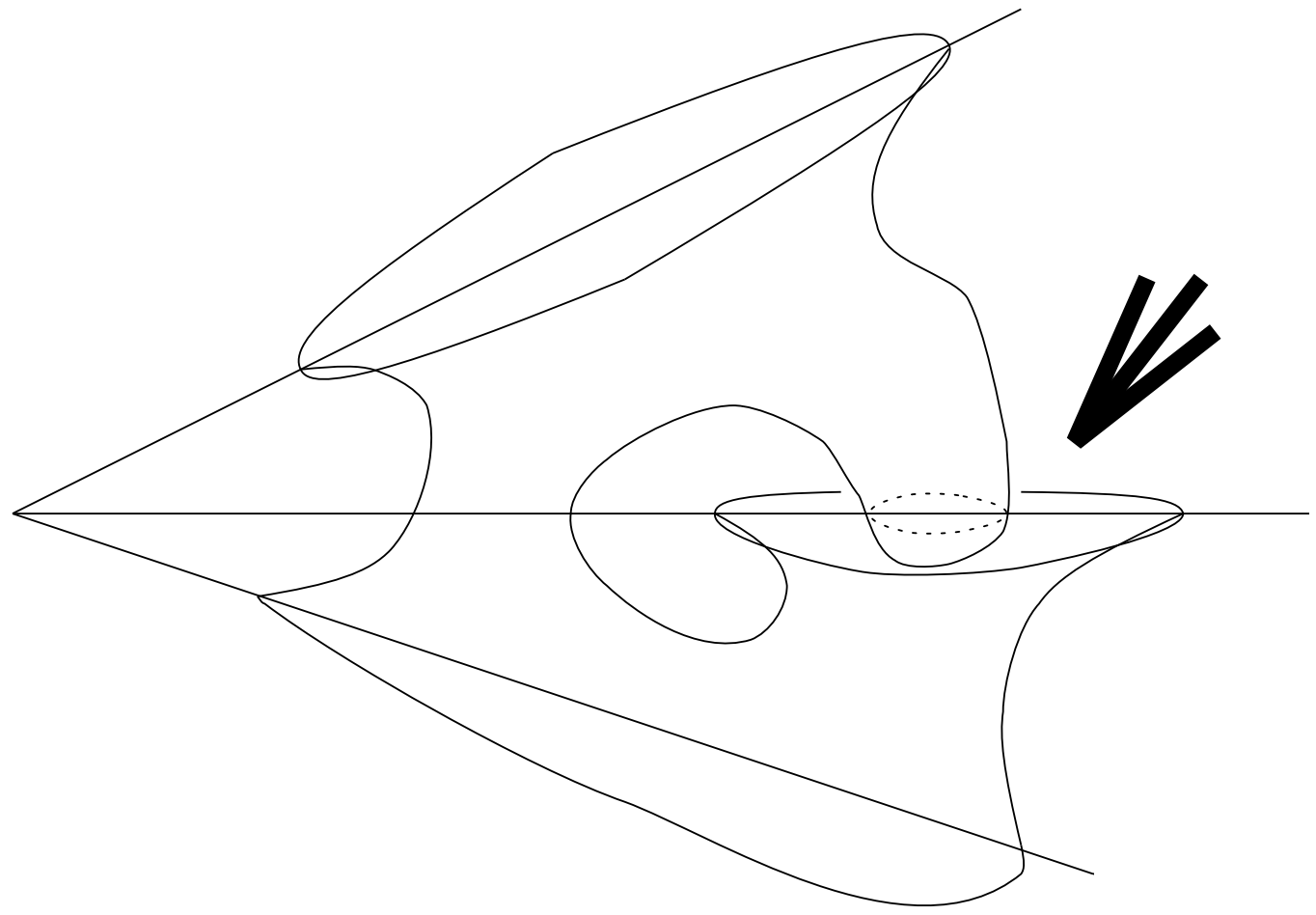

FIGURE 9. A Surface Ruled Out by Condition A1

A2. About each $B_j$, $\mathcal{D}$ lies locally to one side of $\Pi_j$. That is, for each $B_j$, there is a neighborhood $V_j$ of $B_j$ in $\mathbf{R}^3$ and a half space $H_j$ bounded by the plane containing $\Pi_j$ such that $\mathcal{D} \cap V_j \cap H_j = \emptyset$ (see Figure 10).

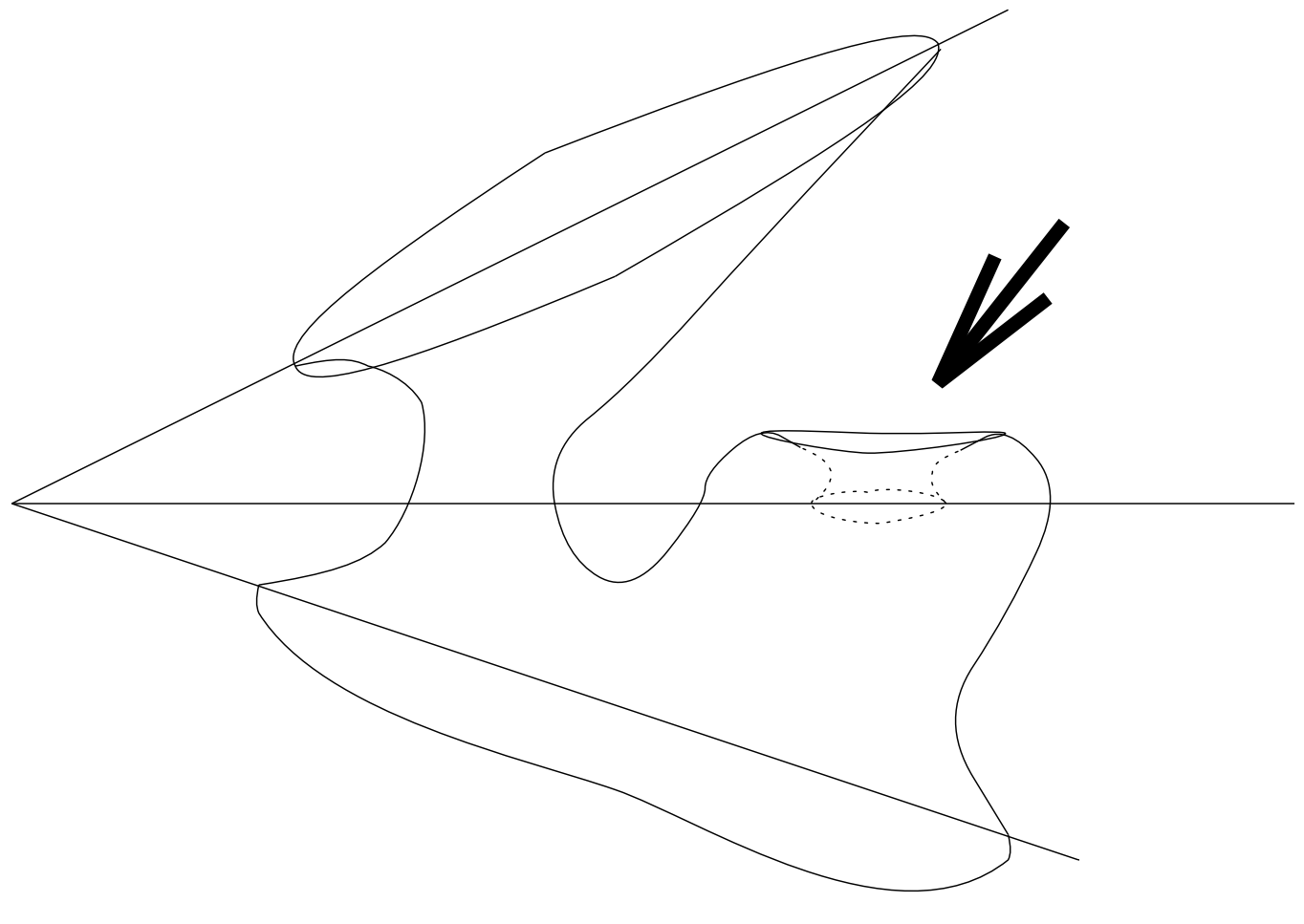

FIGURE 10. A Surface Ruled Out by Condition A2

Finally we require

A3. (contact angle condition) If $N'_j$ is the unit normal to $B_j$ pointing into $\mathcal{D}$ and $N$ is the unit normal to $\mathcal{S}$ pointing into $\mathcal{D}$, then $N'_j \cdot N = c_j$ is a constant along $\beta_j$.

NOTE: When no confusion should arise, we will suppress the index $j$ in regard to $\partial\mathcal{S}$. Notice that with our choice of normals, the constant $c = -\cos\gamma$.
REMARK: Strictly speaking, the wetted regions may be allowed to occupy any part of the plane containing $\Pi$ as long as at least one line through the origin of reflection intersects $\mathcal{S}$ in some set other than a connected segment. It is somewhat natural to ask if this condition can be relaxed as well, but the present method does not seem to apply.

For such a surface $\mathcal{S}$, let us consider carefully the reflection procedure outlined in §2. Let the origin $0 \in \mathbf{R}^3$ be any fixed point on the vertex of the wedge. Recall that we begin with a sphere $B_\rho(0) \supseteq \overline{\mathcal{D}}$ and decrease its radius to a first value $\rho_0$ where $\mathcal{S} \cap S_{\rho_0}(0) \neq \emptyset$. For $\rho < \rho_0$ we wish to apply the map (2.1) to obtain a reflected "comparison" surface $\hat{\mathcal{S}}$.
PROOF of Theorems 4 and 6: It is natural to think of NT1-NT4 as conditions on $\rho$ and on the sets

$$\begin{aligned} NT_1(\rho) &= \{X \in \text{ int } \mathcal{S} : |X| > \rho\}, \\ NT_2(\rho) &= \{X \in \partial\mathcal{S} : |X| > \rho\}, \\ NT_3(\rho) &= \{X \in \text{ int } \mathcal{S} : |X| = \rho\}, \text{ and} \\ NT_4(\rho) &= \{X \in \partial\mathcal{S} : |X| = \rho\} \end{aligned}$$

in particular. Notice that these sets give a disjoint decomposition of $\mathcal{S}^-(\rho)$ for any $\rho$. Let us also add the conditions $X \cdot N < 0$ and $X \cdot n < 0$ to conditions NT1 and NT2 respectively.[5] The proof then reduces, by compactness, to showing the following local result.

**Theorem 17.** *If $P$ satisfies one of* NT1-NT4($\rho$), *then there is some neighborhood $\mathcal{U}$ of $P$ in $\mathcal{S}$ and some $\epsilon > 0$ such that each $X \in \mathcal{U}$ satisfies one of* NT1-NT4($\rho - \epsilon$).[6]

We give the details of the reduction in two steps.
STEP 1. We claim that for some $\epsilon > 0$, every $X$ in $\mathcal{S}^-(\rho_0 - \epsilon)$ satisfies one of NT1-NT4($\rho_0 - \epsilon$).

$\mathcal{S}^-(\rho_0) = \text{NT}_3(\rho_0) \cup \text{NT}_4(\rho_0)$. If $X \in \text{NT}_3(\rho_0)$, then $\mathcal{S}$ is tangent to $S_{\rho_0}$ at $X$. Consequently, $X \cdot N = -|X| < 0$. Therefore, NT3($\rho_0$) is satisfied. If $X \in \text{NT}_4(\rho_0)$, then $\partial\mathcal{S}$ is tangent to $S_{\rho_0} \cap \Pi$ at $X$. Consequently, $X \cdot n = -|X| < 0$. Therefore, NT4($\rho_0$) is satisfied.

Hence, the compact set $\mathcal{S}^-(\rho_0)$ is covered by neighborhoods $\mathcal{U} = \mathcal{U}(X)$ from Theorem 17. Let $\mathcal{U}_1, \dots, \mathcal{U}_k$ be finitely many of these neighborhoods that cover $\mathcal{S}^-(\rho_0)$ and let $\epsilon_1, \dots, \epsilon_k$ be the associated $\epsilon$'s from Theorem 17.

$\mathcal{S} \setminus \cup\mathcal{U}_i$ is a compact set which is disjoint from $S_{\rho_0}$. Let $d$ be the distance between these two sets. Take $\epsilon = \min\{d/2, \epsilon_1, \dots, \epsilon_k\}$.

[5] If the theorem is true, then these conditions evidently follow from NT3 and NT4.
[6] Notice that NT$k(\rho - \epsilon)$ implies NT$k(\rho')$ for $\rho - \epsilon \leq \rho' \leq \rho_0$.

If $X \in \mathcal{S}^-(\rho_0 - \epsilon)$, then $X \in \mathcal{U}_i$ for some $\mathcal{U}_i$. Therefore, $X$ satisfies one of NT1-NT4($\rho_0 - \epsilon_i$). By the footnote to Theorem 17, $X$ satisfies one of NT1-NT4($\rho_0 - \epsilon$).

STEP 2. Let $\rho_1 = \inf\{\rho$: every $X \in \mathcal{S}^-(\rho)$ satisfies one of NT1-NT4($\rho$)$\}$. By STEP 1 $\rho_1 < \rho_0$. By examining points on $\partial\mathcal{S}$ with respect to NT2 we see that $0 < \rho_1$.

We claim that one of NT1-NT4($\rho_1$) must not be satisfied by some point $X \in \mathcal{S}^-(\rho_1)$.

Otherwise, we can take finitely many open sets $\mathcal{U}_1, \dots, \mathcal{U}_k$ from Theorem 17 as before. Letting $\epsilon = \min\{d/2, \epsilon_1, \dots, \epsilon_k\}$ where $d = \text{dist}(\mathcal{S}\setminus \cup \mathcal{U}_i, S_{\rho_1})$, we see that each $X \in \mathcal{S}^-(\rho_1 - \epsilon)$ satisfies one of NT1-NT4($\rho_1 - \epsilon$) which, in this case, is a contradiction.

We may, therefore, consider each of the possible failures in turn.

If NT1($\rho_1$) is not satisfied at $P \in \text{NT}_1(\rho_1)$, then it is still the case, for some neighborhood $\mathcal{U}$ of $P$ in $\mathcal{S}$, that each $X \in \mathcal{U}$ satisfies $(1-\delta)X \in \mathcal{D}$ for $\delta \in (0, 1 - \rho_1^2/|X|^2)$. Therefore, $\hat{X} = (\rho_1^2/|X|^2)X \in \mathcal{D}$ or $\partial\mathcal{D}$. In the latter case, since A2 rules out the possibility that $\hat{X} \in \bar{B}$, $\hat{X} \in \text{int } \mathcal{S}$. It is evident in particular, since $\mathcal{D}$ is an open set, that $\hat{P} \in \text{int } \mathcal{S}$.

Shifting our frame of reference to orthonormal coordinates $x, z$ with $x \in T_{\hat{P}}\mathcal{S} = T_{\hat{P}}\mathcal{S}^+(\rho_1)$ and $N(\hat{P})$ lying along the positive $z$ direction, $\mathcal{S}$ is given locally as the graph of a function $u = u(x)$. It is well known that the map $\psi : (x, z) \mapsto (x, z - u(x))$ is a local diffeomorphism of neighborhoods $\mathcal{V}_1 \longrightarrow \mathcal{V}_2$ of $\hat{P}$ in $\mathbf{R}^3$. Moreover, we can require that $\mathcal{V}_1 \cap \mathcal{D}$ maps under $\psi$ to the "positive $z$"-half space. For $\mathcal{U}$ small enough, the reflection $\hat{\mathcal{U}}$ of $\mathcal{U}$ will lie in $\mathcal{V} \cap (\mathcal{D} \cup \mathcal{S})$. Thus, $\psi\hat{\mathcal{U}}$ is a neighborhood of a smooth surface in the "positive $z$"-half space which passes through $\hat{P} = \psi(\hat{P}) = (0, 0)$. Clearly, then, $\psi\hat{\mathcal{U}}$ is tangent to $T_{\hat{P}}\mathcal{S}$ at $\hat{P}$. Consequently, $\hat{\mathcal{U}} \subset \hat{\mathcal{S}}$ is tangent to $\mathcal{S}^+$ at $\hat{P}$. Moreover, we note the following important fact.

**T1a.** If $\mathcal{S}^+$ and $\hat{\mathcal{S}}$ are expressed locally as graphs of functions $u$ and $\hat{u}$ over a neighborhood $\eta$ of their common tangent plane $T_{\hat{P}}\mathcal{S}$ (oriented as above), then $u \leq \hat{u}$ near $\hat{P}$.

If NT2($\rho_1$) is not satisfied at $P \in \text{NT}_2(\rho_1)$, then it is still the case, for some neighborhood $\mathcal{U}$ of $P$ in $\mathcal{S}$, that each $X \in \mathcal{U}$ satisfies $(1-\delta)X \in \mathcal{D} \cup B$ for $\delta \in (0, 1 - \rho_1^2/|X|^2)$. Focusing our attention on the points on $\beta$ near $P$, we can show, by an argument similar to that for the T1 case above (except in one less dimension), that the reflection of $\beta$, $\hat{\beta}$, is tangent to $\beta$ at $\hat{P}$. Letting $v$ be the unit vector in $T_P\mathcal{S}$ pointing into $\mathcal{S}$ and perpendicular to $\beta$ at $P$, we see that $v \cdot n = \cos\gamma$, and by the conformality of the reflection, $\hat{v} \in T_{\hat{P}}\mathcal{S}$ satisfies $\hat{v} \cdot \hat{n} = \cos\gamma$ and is perpendicular to $\hat{\beta}$. Thus, $\hat{\mathcal{S}}$ is tangent to $\mathcal{S}^+$ at $\hat{P}$ as before. The map $\psi$ in this case allows us to define a curve $b = \psi \circ \beta$ near $\hat{P}$, and by choosing the positive $x_2$-axis along $\hat{v}$ we have that $b$ is tangent to the $x_1$-axis. Thus, $b$ is locally the graph of a function $x_2 = b(x_1)$, and $\mathcal{S}$ is given locally as a graph over $\{x : x_2 \geq b(x_1)\}$ in $T_{\hat{P}}\mathcal{S}$. A similar discussion yields a curve $\hat{b} = \psi \circ \hat{\beta}$ such that $\hat{\mathcal{S}}$ is given locally as a graph over $\{x : x_2 \geq \hat{b}(x_1)\}$ in $T_{\hat{P}}\mathcal{S}$. Both surfaces are represented as graphs over $\eta = \{x : x_2 \geq b_0(x_1) \equiv \max\ \{b(x_1), \hat{b}(x_1)\}\}$,[7] and by the same use of $\psi$ as in the T1 case we see that $\hat{\mathcal{S}}$ lies locally above $\mathcal{S}^+$. In summary,

**T2a.** If $\mathcal{S}^+$ and $\hat{\mathcal{S}}$ are expressed locally as graphs of functions $u$ and $\hat{u}$ over $\eta = \{x : x_2 \geq b_0(x_1)\}$ in $T_{\hat{P}}\mathcal{S}$ as above, then $u \leq \hat{u}$.

The arguments proceed similarly in the T3 and T4 cases except that $b$ is determined by the projection of $S_{\rho_1} \cap \mathcal{S}$ onto $T_P\mathcal{S} = T_{\hat{P}}\mathcal{S}$ in the T3 case, and two curves are necessary in the T4 case. The following statements summarize the situation.

If NT3($\rho_1$) fails at $X \in \mathrm{NT}_3(\rho_1)$, then T3($\rho_1$) holds for $\hat{X} = X$, and

**T3a.** For an appropriate choice of $x, z$-coordinates at $X = \hat{X}$ with $N(X)$ along the positive $z$-axis, $\mathcal{S}^+$ and $\hat{\mathcal{S}}$ are expressed locally as the graphs of functions $u$ and $\hat{u}$ over $\eta = \{x : x_2 \geq b_1(x_1)\} \subset T_X\mathcal{S}$ where $\{b_1\}$ is the projection $\psi(S_{\rho_1} \cap \mathcal{S})$ near $X$. Moreover, $u \leq \hat{u}$ on $\eta$.

If NT4($\rho_1$) fails at $X \in \mathrm{NT}_4(\rho_1)$, then T4($\rho_1$) holds for $\hat{X} = X$, and

**T4a.** Let $b_0(x_1) = \max\{b(x_1), \hat{b}(x_1)\}$ as in T2a and let $\{b_1\}$ be the projection near $X$ of $S_{\rho_1} \cap \mathcal{S}$ onto $T_X\mathcal{S}$. $b_1$ can be expressed as a graph[8] $x_1 = b_1(x_2)$ near $X$, and $b_0'(0) = 0 = b_1'(0)$. $\mathcal{S}^+$ and $\hat{\mathcal{S}}$ can be expressed locally as graphs of functions $u$ and $\hat{u}$ on $\eta = \{x : \pm x_1 \geq \pm b_1(x_2), x_2 \geq b_0(x_1)\}$,[9] and $u \leq \hat{u}$.

This completes the reduction to Theorem 17.

Our next reduction is to the following explanation of the local behavior of the inward normal $N$.

**Theorem 18.** *If $P \in \mathrm{int}\ \mathcal{S}$ and $P \cdot N(P) < 0$, then there is some neighborhood $\mathcal{U}$ of $P$ in $\mathcal{S}$ such that $X \cdot N < 0$ for every $X \in \mathcal{U}$.*

*If $P \in \partial\mathcal{S}$ and $P \cdot n(P) < 0$, then there is some neighborhood $\mathcal{U}$ of $P \in \mathcal{S}$ such that $X \cdot n < 0$ for every $X \in \mathcal{U} \cap \partial\mathcal{S}$ and $X \cdot N < 0$ for every $X \in \mathcal{U} \cap \mathrm{int}\ \mathcal{S}$.*

As the properties described in this theorem will be used repeatedly, we give some notation which saves space. The normal $N$ at a point $X \in \mathrm{int}\ \mathcal{S}$ (or the point $X$ itself) is said to be *central* if $X \cdot N < 0$. Similarly, $n = n(X)$ for $X \in \partial\mathcal{S}$ is *central* if $X \cdot n < 0$. We will also need the following properties of $N$ and $n$.

**Lemma 19.** *Let $v$ be a unit vector. If $X \in \mathrm{int}\ \mathcal{S}$ and $v \cdot N < 0$, then there is some $\delta_0 > 0$ such that $X - \delta v \in \mathcal{D}$ and $X + \delta v \in \mathcal{D}^c$ for every $\delta \in (0, \delta_0]$.*

---

[7] $b_0(x_1) \equiv \max\ \{b(x_1), \hat{b}(x_1)\}$ is really just one of the curves being compared because $\hat{b}$ is locally in the closure of $B$ and $b = \partial B$. In fact, if $\gamma \leq \pi/2$, then $b_0 = \hat{b}$, and if $\gamma \geq \pi/2$, then $b_0 = b$.

[8] This is explained more fully below.

[9] The + sign is taken if $X \cdot (1, 0) = -|X| < 0$. The minus sign is taken if $X \cdot (1, 0) = |X| > 0$.

PROOF: Introduce new coordinates $(x_1, x_2, x_3)$ with origin at $X$ such that $T_X\mathcal{S}$ is the $x_1, x_2$-plane and $N$ is a positive orientation. Let $v = (v_1, v_2, v_3)$ in these new coordinates. The condition $v \cdot N < 0$ evidently becomes, in the new coordinates, $v_3 < 0$.

As is well known, $\mathcal{S}$ is locally given as the graph of a function $u = u(x_1, x_2)$ over $T_X\mathcal{S}$, and the map

$$(x_1, x_2, x_3) \longmapsto (x_1, x_2, x_3 - u(x_1, x_2)) \tag{4}$$

is a local diffeomorphism of neighborhoods of $(x_1, x_2, x_3) = (0, 0, 0)$ in $\mathbf{R}^3$ such that $\mathcal{S}$ maps into the $x_1, x_2$-plane, $x_3 > 0$ corresponds locally to a portion of $\mathcal{D}$ and $x_3 < 0$ corresponds locally to a portion of $\mathcal{D}^c$.

Now, $x - \delta v$ is given in our new coordinates by $-\delta v$, and the image of this point under (4) has $x_3$-component

$$x_3 = -\delta v_3 - u(-\delta v_1, -\delta v_2) = -\delta v_3 + \circ(\delta). \tag{5}$$

For $\delta > 0$ small enough, this point satisfies $x_3 > 0$ and, hence, corresponds to the point $x - \delta v \in \mathcal{D}$.

**Corollary 20.** *If $X \in \text{int } \mathcal{S}$ and $v$ is a unit vector such that for some sequence of $\delta_j \searrow 0$, $X + \delta_j v \in \mathcal{D}$ or $X - \delta_j v \in \mathcal{D}^c$, then $v \cdot N \geq 0$.*

PROOF: Assume the statement is false. Lemma 19 gives an immediate contradiction.

Similarly we have

**Lemma 21.** *Let $v$ be a unit vector in $\Pi$. If $X \in \partial\mathcal{S}$ and $v \cdot n < 0$, then there is some $\delta_0 > 0$ such that $X - \delta v \in B$ and $X + \delta v \in B^c$ for every $\delta \in (0, \delta_0]$.*

**Corollary 22.** *If $X \in \partial\mathcal{S}$ and $v$ is a unit vector such that for some sequence of $\delta_j \searrow 0$, $X + \delta_j v \in B$ or $X - \delta_j v \in \Pi \backslash B$, then $v \cdot n \geq 0$.*

PROOF of Theorem 17: As discussed above, there is no generality lost in adding the condition $X \cdot N < 0$ to NT1 and the condition $X \cdot n < 0$ to NT2.

Say $P$ satisfies one of NT1-NT4($\rho$). Let $\mathcal{U} = \mathcal{U}(P)$ be the neighborhood of $P$ in $\mathcal{S}$ given by Theorem 18. For each $X \in \mathcal{U} \cap \text{int } \mathcal{S}$ there is some $\delta_0 > 0$ such $(1-\delta)X \in \mathcal{D}$ for $\delta \in (0, \delta_0)$. This follows from Lemma 19. By taking the supremum over all such $\delta_0$, we may assume that $\delta_0 = \delta_0(X)$ is taken as large as possible. Similarly, for $X \in \mathcal{U} \cap \partial\mathcal{S}$ we can take the largest possible $\delta_0 = \delta_0(X)$ such that $(1 - \delta)X \in B$ for $\delta \in (0, \delta_0)$. That $\delta_0 > 0$ follows from Lemma 21. Notice that in either case $Y = (1 - \delta_0)X \in \mathcal{S}$. Moreover, by Corollary 20 and Corollary 22, if $X \in \text{int } \mathcal{S}$, then $Y \cdot N(Y) \geq 0$ and if $X \in \partial\mathcal{S}$, then $Y \cdot n(Y) \geq 0$.

For $\eta > 0$ let $\mathcal{U}_\eta = B_\eta(P) \cap \mathcal{S}$. Notice that for $\eta$ small enough $\mathcal{U}_\eta \subset \mathcal{U}$. We claim that if

$$\lim_{\eta \searrow 0} \inf_{X \in \mathcal{U}_\eta} \delta_0(X) > 1 - \rho^2/|P|^2 \tag{6}$$

then for $\eta$ small enough there is some $\epsilon > 0$ such that each $X \in \mathcal{U}_\eta$ satisfies one of NT1-NT4($\rho - \epsilon$). In fact, one of NT3($\rho - \epsilon$) or NT4($\rho - \epsilon$) clearly follows for points with $|X| = \rho - \epsilon$ by Lemma 18 as long as $\eta$ is small enough. On the other hand, since $\sup_{X \in \mathcal{U}_\eta} |X - P| \to 0$ as $\eta \to 0$, if $\eta$ is small enough, then according to (6)

$$\inf_{X \in \mathcal{U}_\eta} \delta_0(X) > \sup_{X \in \mathcal{U}_\eta} (1 - \rho^2/|X|^2).$$

Carrying this one step further, since

$$\lim_{\epsilon \to 0} \sup_{X \in \mathcal{U}_\eta} \left(1 - \frac{(\rho - \epsilon)^2}{|X|^2}\right) = \sup_{X \in \mathcal{U}_\eta} \left(1 - \frac{\rho^2}{|X|^2}\right),$$

we see that for some $\epsilon > 0$

$$\inf_{X \in \mathcal{U}_\eta} \delta_0(X) > \sup_{X \in \mathcal{U}_\eta} \left(1 - \frac{(\rho - \epsilon)^2}{|X|^2}\right).$$

Thus, for any particular $X \in \mathcal{U}_\eta$, $\delta_0(X) > 1 - (\rho - \epsilon)^2/|X|^2$. Our claim then follows from the definition of $\delta_0$.

It remains to verify (6). If (6) were false, then we obtain a sequence $X_j \to P$ with $\delta_{00} = \lim_{j \to \infty} \delta_0(X_j) \le 1 - \rho^2/|P|^2$. Recall that $Y_j \equiv (1 - \delta_0(X_j))X_j \in \mathcal{S}$ and has a non-central normal. Thus, $Y = \lim_{j \to \infty} Y_j = (1 - \delta_{00})P \in \mathcal{S}$. If $\delta_{00} \in (0, 1 - \rho^2/|P|^2]$, then $|P| > \rho$ and this contradicts the fact that $P$ satisfies one of NT1-NT4($\rho$). The only other possibility is that $\delta_{00} = 0$. In that case $Y = P$ and $Y \cdot N(Y) \ge 0$ (or $Y \cdot n(Y) \ge 0$ if $P \in \partial\mathcal{S}$). This contradicts the fact that $P$ has a central normal.

This completes the reduction to Theorem 18.

PROOF of Theorem 18: If $P \in \text{int}\ \mathcal{S}$ and is central, then it follows by the continuity of $N$ that points nearby $P$ in $\mathcal{S}$ are central. Similarly, if $P \in \partial\mathcal{S}$ and is central, then the points in $\partial\mathcal{S}$ near $P$ are central. The following observation gives an easy extension to a full neighborhood of $\mathcal{S}$ when $0 < \gamma < \pi$.

**Lemma 23.** *If $P \in \partial\mathcal{S}$ is central and $0 < \gamma < \pi$, then $P \cdot N(P) < 0$.*

PROOF: Since $n$ and $N'$ span the plane orthogonal to $\beta$,

$$P \cdot N(P) = (N(P) \cdot n(P))(P \cdot n(P)).$$

Therefore, Lemma 23 follows from

**Lemma 24.** *$N \cdot n \ge 0$ on $\partial\mathcal{S}$ with equality if and only if $\gamma = 0$ or $\pi$. If $\gamma = 0$, then $N = -N'$. If $\gamma = \pi$, then $N = N'$.*

PROOF: This result follows essentially from the facts that $\bar{\mathcal{D}} \cap V = \emptyset$ and $B \subset \bar{\mathcal{D}}$; see Figure 7. For convenience, let us express vectors in terms of the orthonormal basis $\{n, N'\}$. $N$ lies in the plane normal to $\{\beta\}$ and, consequently, $\mathcal{S}$ intersects that plane in some curve $C = C(s) \in \text{int}\ \mathcal{S}$ near $X$. If $N \cdot N' > 0$, but $N \cdot n < 0$, then $C$ enters the first quadrant (because it cannot enter $V$). By considering $C - \delta N$ which for some $\delta_0 = \delta_0(s)$ is in $B$, we obtain a sequence of points $X_j \in \mathcal{S}$ which converge to $X$ and satisfy $X_j - \delta N \in \mathcal{D}$ for $\delta \in (0, \delta_1)$ and $X_j - \delta_1 N = C(s) - \delta_0(s) N \in B$ for some $s$. By Corollary 20, $N(X_j) \cdot N \leq 0$. Since $N(X_j) \to N$, this is a contradiction.

The other cases follow similarly.

Thus, the proof of Theorem 18 will be complete if we can give an analysis for the extreme angles 0 and $\pi$.

**Theorem 25.** *If $P \in \partial\mathcal{S}$ and $P \cdot n(P) < 0$ and $\gamma = 0$ or $\pi$ at $P$, then there is a neighborhood $\mathcal{U}$ of $P$ in $\mathcal{S}$ such that $X \cdot N(X) < 0$, i.e., $X$ is central, for every $X \in \mathcal{U} \cap \operatorname{int} \mathcal{S}$.*

PROOF: Assume $\gamma = 0$. The argument for $\gamma = \pi$ is similar.

Let $\beta$ be parameterized near $P$ by its arclength $s$. For an appropriate orientation of $\beta$, the vectors $\dot{\beta}(0)$, $n(0)$ and $N'$ constitute a positively oriented orthonormal frame at $P$. Let us express points in coordinates $(x, z)$ with respect to this frame where $x = (x_1, x_2) \in \mathbf{R}^2$ and $z \in \mathbf{R}^1$. Since $T_P\mathcal{S} = \Pi$ and $P \cdot n(P) < 0$, the center of reflection $0 \in l$ is given in these new coordinates by $(\xi, 0)$ for some $\xi \in \mathbf{R}^2$ with $\xi_2 > 0$.

It will be convenient for us to use square balls $\mathrm{B}_r(P) \equiv \{x : -r < x_1, x_2 < r\}$ as well as the standard balls $B_r(P)$.

We claim that there exist smooth real valued functions, $b$ defined in a neighborhood of $0 \in \mathbf{R}$ and $u$ defined in a neighborhood of $0 \in \mathbf{R}^2$, which satisfy $b'(0) = 0 = b(0)$ and $u(x_1, b(x_1)) = 0 = \nabla u(x_1, b(x_1))$ and for which, if $r$ is small enough, the following hold.

- **(i):** $\partial\mathcal{S} \cap \mathrm{B}_r(P) = \{(x_1, b(x_1)) : |x_1| < r\}$.
- **(ii):** $B \cap \mathrm{B}_r(P) = E_r \equiv \{x \in \mathrm{B}_r(P) : x_2 > b(x_1)\}$.
- **(iii):** $(\Pi \backslash \bar{\mathcal{D}}) \cap \mathrm{B}_r(P) = F_r \equiv \{x \in \mathrm{B}_r(P) : x_2 < b(x_1)\}$.
- **(iv):** $\mathcal{S} \cap B_r(P) \subset \{(x, u(x)) : x \in \bar{E}_r\} \subset \mathcal{S}$ and $u(x) \geq 0$ for $x \in \bar{E}_r$.
- **(v):** $\mu(x_1) \equiv (\xi - \beta(s)) \cdot n(s) > (1/2)\xi_2 > 0$ for $|x_1| \leq r$ where $s = s(x_1) = \int_0^{x_1} \sqrt{1 + b'^2}$.

NOTES: (i)-(iii) follow from the zero contact angle; see Figure 11. (iv), that $\mathcal{S}$ is given locally as a graph, also follows from zero contact angle. That $u \geq 0$ follows from the adherence condition A2. (v) follows from continuity since $\mu(0) = \xi_2$.

Let us consider such a neighborhood determined by $r < |\xi|/2$ and calculate the normal curvature of a curve in an arbitrary direction $v \in \Pi$ at $\beta(s)$. A suitable curve for this calculation is given by

$$\alpha(t) = (\beta(s) + tv, u(\beta(s) + tv)). \tag{7}$$

It will be convenient to express vectors in terms of the adapted frame $\{\dot\beta(s), n(s)\}$. In particular, we write

$$\begin{aligned}\nabla u &= (\nabla u\cdot\dot\beta)\dot\beta + (\nabla u\cdot n)n\\ &= u_1^s\dot\beta + u_2^s n;\\ v &= v_1^s\dot\beta + v_2^s n.\end{aligned}$$

With this notation we see that $\alpha' = (v, u_1^s v_1^s + u_2^s v_2^s)$. Therefore, since

$$\dot\alpha = \frac{\alpha'}{|\alpha'|}$$

where "$\dot{}$" denotes differentiation with respect to arclength, we will need to take the second directional derivatives

$$\begin{aligned}\nabla u_1^s &= (\nabla u_1^s\cdot\dot\beta)\dot\beta + (\nabla u_1^s\cdot n)n\\ &= u_{11}^s\dot\beta + u_{12}^s n;\\ \nabla u_2^s &= u_{21}^s\dot\beta + u_{22}^s n.\end{aligned}$$

Using the fact that $\nabla u \equiv 0$ along $\beta$, we see that

$$\ddot\alpha(0) = (0, u_{11}^s(v_1^s)^2 + 2u_{12}^s v_1^s v_2^s + u_{22}^s(v_2^s)^2).$$

Thus, the normal curvature $k_v$ is given by

$$k_v = -(u_{11}^s(v_1^s)^2 + 2u_{12}^s v_1^s v_2^s + u_{22}^s(v_2^s)^2). \tag{8}$$

Taking $v = \dot\beta$ we have

$$k_{\dot\beta} = -u_{11}^s = -\frac{\partial^2 u}{\partial\dot\beta^2}(\beta(s)).$$

Clearly, the normal curvature of $\beta$ is 0. Thus, $u_{11}^s \equiv 0$.

On the other hand, since $\mathcal{S}$ meets the plane $\Pi$ in a constant contact angle, we may apply Joachimstahl's theorem [3, pg. 152] to conclude that $\beta$ is a curvature line for $\mathcal{S}$. Therefore, taking $v = n$, the other principal direction,

$$k_n = -u_{22}^s = 2H.$$

If we assume that $H \geq 0$, then we obtain a contradiction from the Hopf boundary point comparison theorem (Theorem 28). Therefore, $u_{22}^s \equiv -2H > 0$ along $\beta$.

Finally,

$$u_2^s = \frac{\partial u}{\partial n}(\beta(s)) \equiv 0.$$

This allows us to calculate $u_{12}^s = u_{21}^s$ "along $\beta$." That is

$$\begin{aligned} 0 &= \lim_{h\to 0}\frac{1}{h}\{u_2^{s+h}-u_2^s\} \\ &= \lim_{h\to 0}\frac{1}{h}\{\frac{\partial u}{\partial n(s)}(\beta(s+h))-\frac{\partial u}{\partial n(s)}(\beta(s))\} \\ &= \lim_{h\to 0}\nabla u_2^s(\beta(s^*))\cdot\dot{\beta}(s^*) \end{aligned}$$

for some $s^* \in (s, s+h)$. Evaluating the limit on the right we obtain

$$u_{12}^s = 0.$$

Thus, (8) becomes

$$k_v = 2H(v_2^s)^2.$$

Let us take a (perhaps) smaller $r < 1$ satisfying two additional conditions (aside from (i)-(v)). The first of these additional requirements is the simple continuity condition that for $|x| < r$,

$$\frac{1}{|\xi_2 - x_2|} < \frac{2}{\xi_2} < \infty. \tag{9}$$

The second condition arises from Taylor's formula for $\alpha'$. Since $\alpha'(0) = (v, 0)$ and $\alpha''(0) = (0, u_{22}^s(v_2^s)^2)$, we can express $\alpha'$ as

$$\alpha' = (v, u_{22}^s(v_2^s)^2 t + (1/2)T(x^*)t^2) \tag{10}$$

where $x^* \in \mathrm{B}_r(P)$ and $T$ is the third derivative of $u$ along $\alpha$. Following previous notation

$$T = \sum_{i,j,k=1}^{2} v_i v_j v_k u_{ijk}^s$$

which can be expressed in terms of the third derivatives of $u$.

Notice that $|T| \le \sum |u_{ijk}^s|$ which is a continuous function. Thus, setting

$$M = \max_{X\in\overline{\mathrm{B}_r(P)}} \sum |u_{ijk}^s| < \infty,$$

we have a bound on $|T|$ independent of $s$ or $v$. We take $r$ such that if $t < r$, then

$$\frac{1}{2}Mt^2 \le |H|\frac{\xi_2^2}{16|\xi|^2}t. \tag{11}$$

Having chosen $r$, we take $r_1 < r/2$ such that if $|x_1| < r_1$, then

$$|b'(x_1)| < \epsilon \equiv \frac{1}{1 + 24\frac{|\xi|}{\xi_2}}. \tag{12}$$

In particular, $|b(x_1)| = |b(x_1) - b(0)| \le \epsilon|x_1|$.

We are interested in the particular choice

$$v = \frac{\xi - \beta}{|\xi - \beta|} \tag{13}$$

and the associated curve on $\mathcal{S}$ given by (7). Such curves exhaust an entire neighborhood of $P$ in $\mathcal{S}$; see Figure 11.

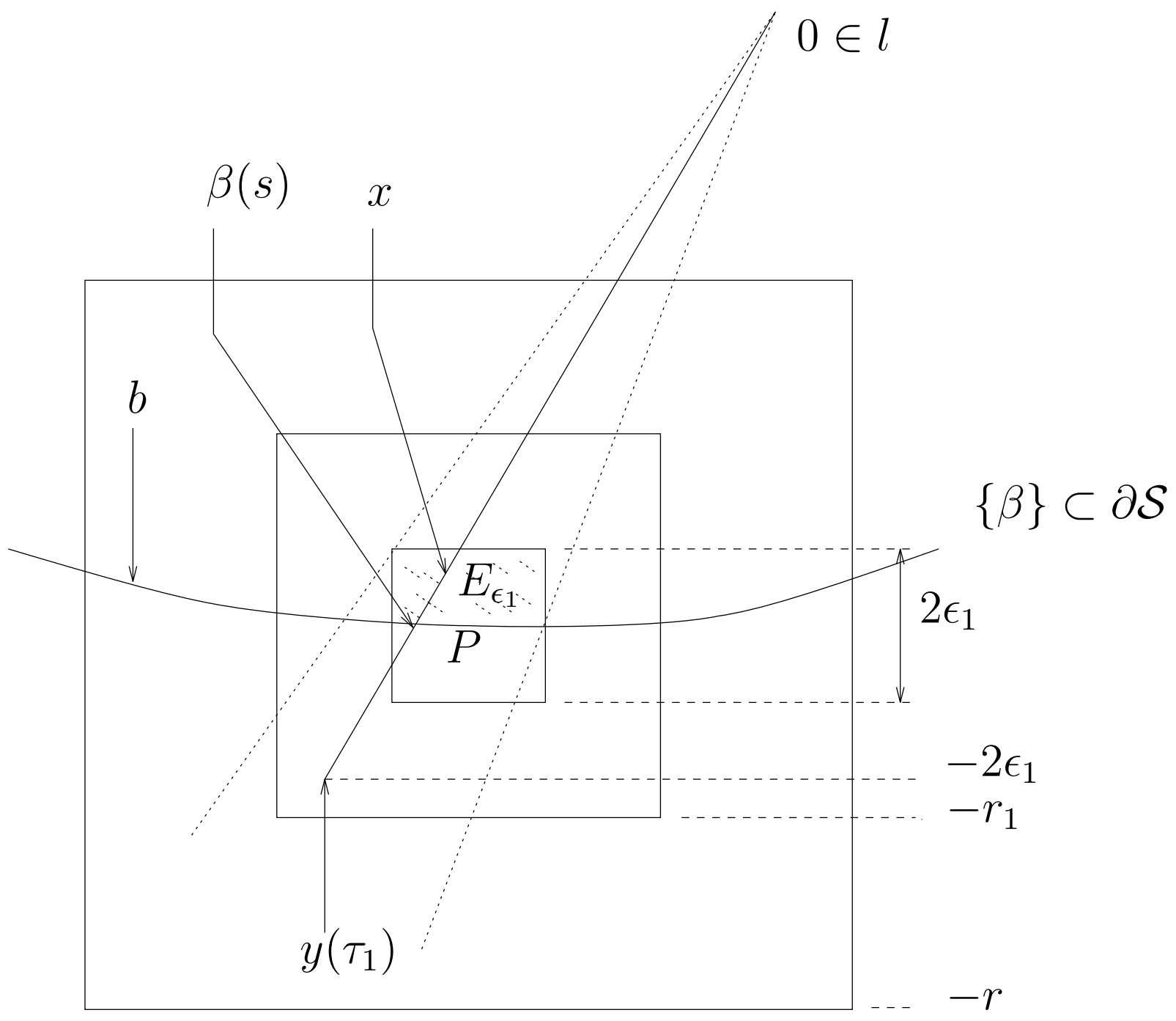

FIGURE 11. A Neighborhood Covered by Rays to the Origin of Reflection

**Claim 26.** *Let*

$$\epsilon_1 = r_1 \min\{1/2, \epsilon\}.$$

*If $x \in E_{\epsilon_1}$, then there is some $\beta(s) \in \mathrm{B}_{r_1}(P)$ such that $(x, u(x))$ lies on the curve $\alpha$ determined by the direction (13). Moreover, there is some $t \in (0, r)$ such that $\alpha(t) = (x, u(x))$.*

PROOF: Let $w = (x - \xi)/|x - \xi|$ and consider the line $y(\tau) = x + \tau w$. $y(0) \in E_{\epsilon_1}$. Therefore,

$$b(y_1(0)) < y_2(0). \tag{14}$$

Setting $\tau_1 = (x_2 - 2\epsilon_1)|\xi - x|/(\xi_2 - x_2)$, we see that $y(\tau_1) \in \mathrm{B}_{r_1}(P)$. In fact,

$$|y_2(\tau_1)| = |-2\epsilon_1| < r_1,$$

and

$$\begin{aligned}
|y_1(\tau_1)| &= \left| x_1 + \frac{(x_2+2\epsilon_1)|x-\xi|}{\xi_2 - x_2}(x_1-\xi_1) \right| \\
&< |x_1| + 8\frac{|\xi|}{\xi_2}3\epsilon_1 \\
&< \left(1 + 24\frac{|\xi|}{\xi_2}\right)\epsilon_1 \qquad (15) \\
&\leq r_1.
\end{aligned}$$

Here we used (9) and the definition of $\epsilon_1$. Thus, according to (12) and (15),

$$|b(y_1(\tau_1))| \leq \epsilon |y_1(\tau_1)| < \epsilon_1.$$

On the other hand, $y_2(\tau_1) = -2\epsilon_1$. We conclude that $y(\tau_1) \in F_{r_1}$. In particular,

$$y_2(\tau_1) < b(y_1(\tau_1)). \qquad (16)$$

Comparing (16) and (14), we see that for some $\tau \in (0, \tau_1)$, $y_2(\tau) = b(y_1(\tau))$.

This is the point we want in $\partial\mathcal{S}$; $\beta(s) = y(\tau)$.

It is clear that $-w = v$ where $v$ is given in (13). Furthermore, $\beta(s) - \tau w = x$, so that we take $t = \tau$. Finally, since $x, y(\tau_1) \in \mathrm{B}_{r_1}(P)$ and $r_1 < r/2$, we see that $t < \sqrt{2}r_1 < r$. This completes the proof of the claim.

In order to finish the proof of Theorem 25, consider $X = (x, u(x))$ in a small enough neighborhood of $P$ in $\mathcal{S}$ so that $x$ is in $E_{\epsilon_1}$. If $N$ is the normal at this point, then

$$\begin{aligned}
X \cdot N &= (x - \xi, u(x)) \cdot (\nabla u, -1) \\
&\leq (x-\xi)\cdot \nabla u \\
&= -|x-\xi| v \cdot \nabla u.
\end{aligned}$$

On the other hand, since $\alpha'$ is tangent to $\mathcal{S}$,

$$\begin{aligned}
0 &= N \cdot \alpha' \\
&= v\cdot\nabla u - (u^s_{22}(v^s_2)^2 t + \frac{1}{2}T(x^*)t^2) \qquad \text{by (10)} \\
&\leq v\cdot\nabla u + 2H(v\cdot n)^2 t + \frac{1}{2}Mt^2 \\
&\leq v\cdot\nabla u + 2H\left(\frac{\mu(s)}{|\xi-\beta|}\right)^2 t - H\frac{\xi_2^2}{16|\xi|^2}t \qquad \text{by (11)} \\
&\leq v\cdot\nabla u + \left(2H\left(\frac{\xi_2}{4|\xi|}\right)^2 - H\frac{\xi_2^2}{16|\xi|^2}\right)t \qquad \text{by (v)} \\
&= v\cdot\nabla u + H\frac{\xi_2^2}{16|\xi|^2}t.
\end{aligned}$$

Thus,

$$\begin{aligned} X \cdot N &\leq |x - \xi| H \frac{\xi_2^2}{16|\xi|^2} t \\ &< 0. \end{aligned}$$

Having established Theorem 4 and the occurrence of T1a-T4a in particular, we note that the fundamental mechanism that allows one to conclude the coincidence of reflected and unreflected portions of a surface is a "comparison theorem." Our treatment of spherical reflection requires the use of three such theorems. We now state these results and check carefully the hypotheses of each application. Furthermore, we give a generalization of Serrin's boundary point lemma at a corner [11] that is used to prove the third comparison theorem. Our treatment of this lemma clarifies the regularity of the coefficients assumed by Serrin and by Gidas, Ni and Nirenberg [5] at least for two dimensional domains with corners of angle $\pi/2$.

For a (T1) touching point we apply[10]

**Theorem 27.** *Let $\mathcal{M}$ be the mean curvature operator, i.e., for a smooth real valued function $u = u(x, y)$ defined on an open set in $\mathbf{R}^2$, $\mathcal{M}u = \mathcal{M}u(x, y)$ is the mean curvature of the graph of $u$ at $(x, y, u(x, y))$ with respect to the upward pointing normal. If $u$ and $\hat{u}$ are defined in a neighborhood of $(x_0, y_0)$ and*

$$u \leq \hat{u}, \tag{17}$$

$$u(x_0, y_0) = \hat{u}(x_0, y_0), \tag{18}$$

$$\mathcal{M}u \geq \mathcal{M}\hat{u}, \tag{19}$$

*then $u \equiv \hat{u}$.*

The tangency is set up by condition T1 and described in the terms of our comparison theorem in T1a. (17) and (18) follow from these statements. (19) evidently follows from the boundedness property when it holds.

The other two comparison theorems apply to an operator defined on the closure of an open connected set $\Omega$ in $\mathbf{R}^2$ bounded by two $C^2$ curves $\gamma_j : [0, s_j] \to \mathbf{R}^2$ which intersect in the point $\gamma_j(0) = x_0$ and an arc of a circle centered at $x_0$ and connecting $\gamma_1(s_1)$ to $\gamma_2(s_2)$.

**Theorem 28** (E. Hopf Boundary Point Comparison Theorem)**.** *Assume $\Omega$ is a half neighborhood, i.e., $\gamma_1$ extends smoothly across $x_0$. i.e., $\gamma_1(-\epsilon) = \gamma_2(\epsilon)$ for $\epsilon$ small and positive.*[11] *Let $\mathcal{M}$ be the mean curvature operator as in Theorem 27. If $u$ and $\hat{u}$ are defined in a neighborhood of $x_0$,*

[10] These results are essentially well known; see also [7].

[11] Here we are assuming that $s$ is an arclength parameter.

$$u \leq \hat{u} \text{ on } \Omega, \tag{20}$$
$$u(x_0) = \hat{u}(x_0), \tag{21}$$
$$\mathcal{M}u \geq \mathcal{M}\hat{u} \text{ on } \Omega, \tag{22}$$

*and*

$$\frac{\partial u}{\partial m}(x_0) = \frac{\partial \hat{u}}{\partial m}(x_0)$$

*where $m$ is the inward normal to $\partial\Omega$ at $x_0$, then $u \equiv \hat{u}$ on $\Omega$.*

**Theorem 29** (Serrin's Corner Comparison Theorem)**.** *Assume $\dot{\gamma_1}(0) \cdot \dot{\gamma_2}(0) = 0$. If $u$ and $\hat{u}$ are defined in a neighborhood of $x_0$, (20), (21) and (22) hold, and*

$$\frac{\partial u}{\partial m}(x_0) - \frac{\partial \hat{u}}{\partial m}(x_0) = 0 = \frac{\partial^2 u}{\partial m^2}(x_0) - \frac{\partial^2 \hat{u}}{\partial m^2}(x_0)$$

*where $m$ is any direction which is non-tangential to $\partial\Omega$ and enters $\Omega$ at $x_0$, then $u \equiv \hat{u}$ on $\Omega$.*

(T2) and (T3) touching points are handled by Theorem 28. The additional condition follows simply from the fact that $\mathcal{S}$ and $\hat{\mathcal{S}}$ are tangent at the point of touching.

Theorem 29 applies in the (T4) case. The agreement of the first order directional derivatives follows from the tangency as before, since $u_m(0) = Du(0) \cdot m = 0$. The second derivative condition is obtained by showing that all the second partial derivatives of $\hat{u}$ and $u$ agree at $0 \in T_X\mathcal{S}$.[12] Let us recall the situation from T4 and T4a (see Figure 12). The surface $\mathcal{S}$ is represented as graph $u$ on a neighborhood $\eta = \{x : x_2 \geq b(x_1)\}$ where graph $b$ is the projection $\psi\{\beta\}$ of the boundary curve $\{\beta\}$ onto $T_X\mathcal{S}$.[13] Assume for definiteness, that $\xi = (\rho_1, 0, 0)$ gives the coordinates of $0 \in l$. Let us consider the intersection of graph $u$ with $S_{\rho_1}$, or more precisely, its projection onto $T_X\mathcal{S}$. This is the set

$$\mathcal{I} = \{x : |x - \xi|^2 + (u(x))^2 = \rho_1^2\}.$$

Since $\mathcal{S}$ is transverse to $S_{\rho_1}$ at $0 \in T_X\mathcal{S}$,

$$\frac{\partial}{\partial x_1}\{|x - \xi|^2 + (u(x))^2\} = 2(x_1 - \rho_1) + 2uu_1, \tag{23}$$

and

$$\frac{\partial}{\partial x_2}\{|x - \xi|^2 + (u(x))^2\} = 2x_2 + 2uu_2,$$

we see that $\mathcal{I}$ is a curve given implicitly as the graph of a function $b_1(x_2) = x_1$ with $b_1(0) = 0 = b_1'(0)$.

[12] Any second directional derivative can be expressed in terms of first and second partial derivatives.

[13] Here $\eta$ is a "half neighborhood" as in Theorem 28, not the "corner domain" given in T4a.

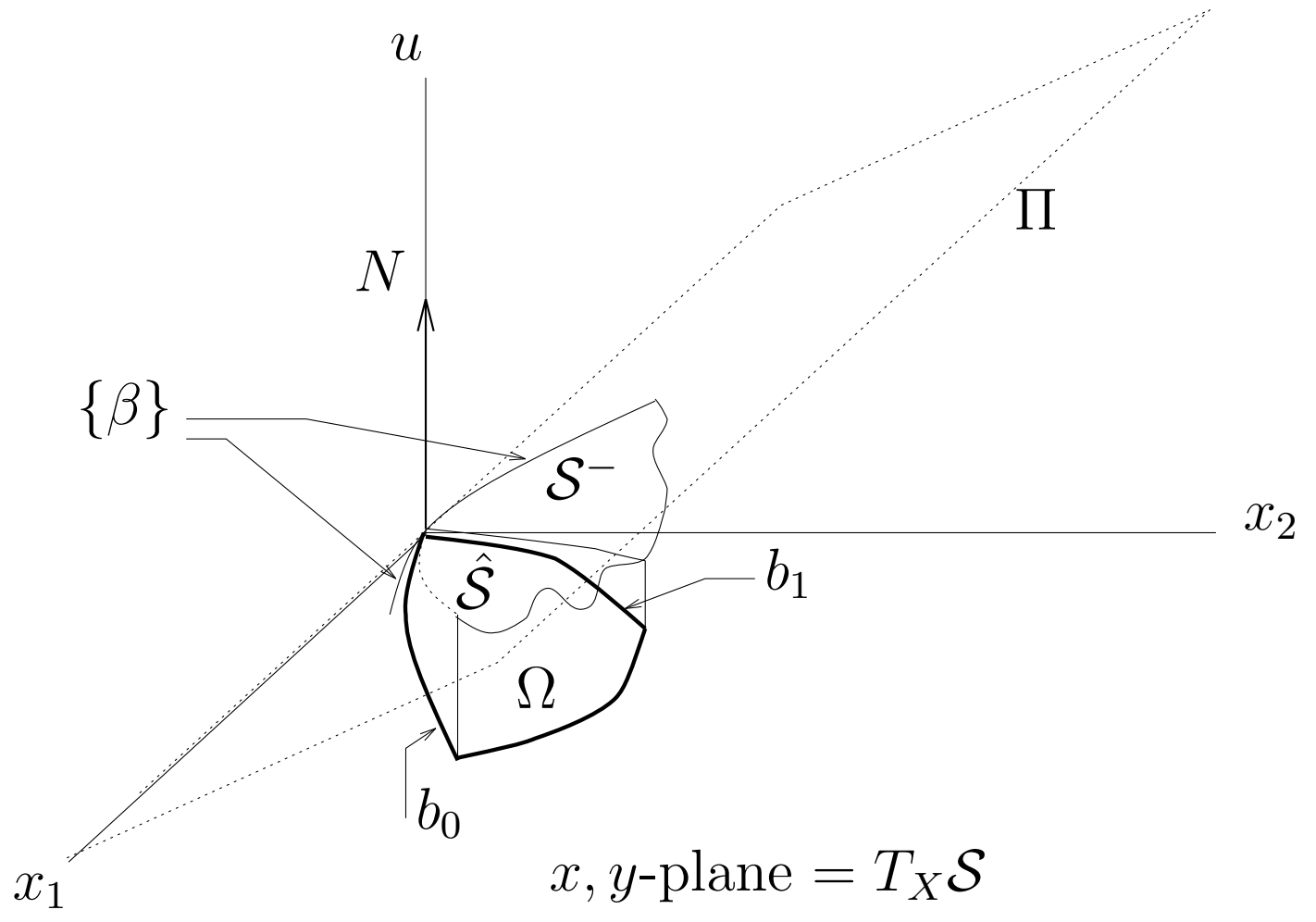

FIGURE 12. A (T4) Touching Point

Furthermore, since (23) is strictly negative near $x = 0$, we see that, for a small enough half neighborhood $\eta$ of 0,

$$G^- \equiv \{(x, u(x)) : x_1 \leq b_1(x_2)\} \subset \mathcal{S}^-(\rho_1)$$

and

$$G^+ \equiv \{(x, u(x)) : x_1 \geq b_1(x_2)\} \subset \mathcal{S}^+(\rho_1).$$

Thus, the comparison of $u$ and $\hat{u}$ makes sense on the corner domain $\Omega = \{x : x_1 \geq b_1(x_2) \text{ and } x_2 \geq b_0(x_1)\}$.

If $R$ is the reflection map restricted to $T_X\mathcal{S}$, then $\hat{u}$ is given on $\Omega$ by $\hat{u} = u \circ R$.[14] $D\hat{u} = Du \circ DR$, and as noted above, this vector vanishes at $x = 0$. One also checks, by direct calculation, that $\hat{u}_{11}(0) = u_{11}(0)$, $\hat{u}_{12}(0) = -u_{12}(0)$, and $\hat{u}_{22}(0) = u_{22}(0)$.[15] On the other hand, we can calculate $u_{12}(0)$ along $\{\beta\}$ as follows.

$$\begin{aligned} \frac{\partial}{\partial x_1}\left[u_2(x_1, b(x_1))\right]\bigg|_{x=0} &= u_{12}(0) + u_{22}(0)b'(0) \\ &= u_{12}(0). \end{aligned}$$

Recall, however, that along this curve the contact angle condition is expressed by

$$-\cos\gamma = N \cdot N' = \frac{1}{\sqrt{1 + |Du|^2}}(-N_2' u_2 + N_3') \tag{24}$$

[14] This replaces the explicit formula $\hat{u}(x) = u(-x_1, x_2)$ that arises in planar reflection.

[15] Notice that $D^2R$ plays no role in these calculations due to multiplication by $Du(0) = 0$.

where $N_2' \geq 0$ with equality holding if and only if $\gamma = 0$ or $\pi$. If $N_2' = 0$, then $\Pi = T_X\mathcal{S}$ and

$$u(x_1, b(x_1)) \equiv 0 = u_2(x_1, b(x_1)) = \frac{\partial}{\partial x_1}\left[u_2(x_1, b(x_1))\right].$$

If $N_2' \neq 0$, then by differentiating (24) we obtain

$$N_2' \frac{\partial}{\partial x_1}\left[u_2(x_1, b(x_1))\right] = \frac{Du \cdot \frac{\partial}{\partial x_1}(Du)}{\sqrt{1+|Du|^2}} N_3'.$$

The right hand side of this equation vanishes with $Du$ at $x = 0$. In either case, $\hat{u}_{12}(0) = -u_{12}(0) = 0$.

For convenience, we state the general principle demonstrated by this argument.

**Lemma 30.** *Let $\mathcal{S}$ be a surface that meets a plane at a constant contact angle along a portion of its boundary near $X \in \partial\mathcal{S}$. If $\mathcal{S}$ is expressed as the graph of a function $u$ on $T_X\mathcal{S}$, then $u_{12}(0) = 0$.*

We conclude this section with some remarks on the derivation of Theorem 29. As outlined in [7] Theorem 29 is derived from the following

**Lemma 31.** *Let $L$ be an elliptic linear operator of the form*

$$L = \sum_{i,j=1}^{2} a_{ij}(x) D_i D_j + \sum_{j=1}^{2} b_j(x) D_j + c(x)$$

*with bounded coefficients satisfying $a_{ij} = a_{ji}$ on the closure of a corner domain $\Omega \in \mathbf{R}^2$ as in Theorem 29. We require, furthermore, that for some constant $K > 0$, the coefficients satisfy*

$$\left|\sum a_{ij}(x)\dot{\gamma_1}^i(0)\dot{\gamma_2}^j(0)\right| \leq K|x - x_0| \tag{25}$$

*for $x \in \Omega$ where $\gamma_j = (\gamma_j^1, \gamma_j^2)$ for $j = 1, 2$.*

*If $u \in C^2(\Omega)$, $u \geq 0$, $Lu \leq 0$, and $u(x_0) = 0 = u_m(x_0) = u_{mm}(x_0)$ for every nontangential direction $m$, then $u \equiv 0$ on $\Omega$.*

REMARK: Serrin [11] proved a similar result under a regularity assumption on the coefficients which was stricter than (25) and under the hypothesis that one of the curves $\gamma_j$ that bound $\Omega$ was a straight line (which is enough for planar reflection). Gidas, Ni, and Nirenberg [5] relaxed the condition on the boundary curve and extended the result to corner domains with angle smaller than $\pi/2$ but assumed that the coefficients were $C^1$. By combining the methods of these two papers, one can obtain the result as stated above. In fact, if $\dot{\gamma_1}(0) = (0, 1)$ and $\dot{\gamma_2}(0) = (1, 0)$, then one obtains the result under the condition $|a_{12}(x)| \leq K(|x_1| + |x_2|)$. (25) can easily be shown a sufficient condition from this case by a change of variables.

John McCuan, Mathematics Department, Stanford University, Stanford, CA 94025
*E-mail address*: johnm@math.stanford.edu